\documentclass[11pt,leqno]{amsart}
\usepackage{amsmath,amsfonts,amssymb,amscd,amsthm,amsbsy,upref,mathrsfs}
\numberwithin{equation}{section}
\newtheorem{thm}{Theorem}[section]
\newtheorem{lem}[thm]{Lemma}
\newtheorem{cor}[thm]{Corollary}
\newtheorem{prop}[thm]{Proposition}
\theoremstyle{definition}

\newtheorem{defn}[thm]{Definition}
\newtheorem{assump}[thm]{Assumption}
\newtheorem{prob}[thm]{Problem}

\newtheorem*{Defn}{Definition}

\newtheorem*{remark}{Remark}

\newcommand{\supp}{\mathrm{\rm supp\,}}
\newcommand{\dist}{\mathrm{dist}}
\newcommand{\ball}{\mathrm{\rm ball\,}}
\newcommand{\ran}{\mathrm{\rm ran\,}}
\newcommand{\rank}{\mathrm{\rm rank\,}}
\newcommand{\age}{\mathrm{\rm age\,}}

\newcommand{\weight}{\mathrm{\rm weight\,}}

\newcommand{\im}{\mathrm{\rm im\,}}
\newcommand{\reals}{\mathbb R}
\newcommand{\nats}{\mathbb N}

\newcommand{\rats}{\mathbb Q}
\newcommand{\half}{{\textstyle\frac12}}
\newcommand{\txtfrac}[2]{{\textstyle\frac{#1}{#2}}}
\newcommand{\XK}{\mathfrak X_{\mathrm K}}
\newcommand{\Xgm}{\mathfrak X_{\mathrm {gm}}}
\newcommand{\Xad}{\mathfrak X_{\mathrm {ad}}}
\newcommand{\BmT}{\mathfrak B_{\mathrm {mT}}}

\begin{document}
\title[The scalar--plus--compact problem] {A hereditarily indecomposable $\mathscr L_\infty$-space that solves the scalar--plus--compact problem}
\author{Spiros A. Argyros and Richard G. Haydon}
\address{Department of Mathematics\\ National Technical University of Athens\\
Greece}

\email{sargyros@math.ntua.gr}
\address{Brasenose College\\ Oxford  OX1 4AJ, U.K.} \email{richard.haydon@brasenose.oxford.ac.uk}
\dedicatory{(23 March 2009)} \subjclass{46B03, 46B26}

\begin{abstract}
We construct a hereditarily indecomposable Banach space with dual space
isomorphic to $\ell_1$.  Every bounded linear operator on this space is
expressible as $\lambda I+K$ with $\lambda$ a scalar and $K$ compact.
\end{abstract}

\maketitle

\section{Introduction}

The question of whether there exists a Banach space $X$ on which
every bounded linear operator is a compact perturbation of a scalar
multiple of the identity has become known as the
``Scalar--plus--Compact Problem''. It is mentioned by Lindenstrauss
as Question 1 in his 1976 list of open problems in Banach space
theory \cite{L}. Lindenstrauss remarks that, by the main theorem of
\cite{ASm} or \cite{Lom}, every operator on a space of this type has
a proper non-trivial invariant subspace. Related questions go
further back: for instance, Thorp \cite{T} asks whether the space of
compact operators $\mathcal K(X;Y)$ can ever be a proper
complemented subspace of $\mathcal L(X;Y)$. On the Gowers--Maurey
space $\Xgm$ \cite{GM}, every operator is a {\em strictly singular}
perturbation of a scalar, and other hereditarily indecomposable (HI)
spaces also have this property. Indeed it seemed for a time that
$\Xgm$ might already solve the scalar--plus--compact problem.
However, after Gowers\cite{G2} had shown that there is a strictly
singular, non-compact operator from a subspace of $\Xgm$ to $\Xgm$,
Androulakis and Schlumprecht \cite{AS} showed that such an operator
can be defined on the whole of $\Xgm$. Gasparis \cite{Ga} has done
the same for the Argyros--Deliyanni space $\Xad$ of \cite{AD}.

In the present paper, we solve the scalar--plus--compact problem by combining
techniques that are familiar from other HI constructions with an additional
ingredient, the Bourgain--Delbaen method for constructing special $\mathscr
L_{\infty}$-spaces \cite{BD}. The initial motivation for combining these two
constructions was to exhibit a hereditarily indecomposable predual of $\ell_1$;
such a space is, in some sense, the extreme example of a known
phenomenon---that the HI property does not pass from a space to its dual
\cite{F,ATol,AF}.  Serendipitously, it turned out that the additional structure
was just what we needed to show that strictly singular operators are compact.
It is interesting, perhaps, to note that the Schur property of $\ell_1$ does
not play a role in our proof and, indeed, we have no general result to say that
an HI predual of $\ell_1$ necessarily has the scalar--plus--compact property.
We use in an essential way the specific structure of the BD construction, which
embeds into our space some very explicit finite-dimensional
$\ell_\infty$-spaces.  As well as the (now) classical machinery of HI
constructions---a space of Schlumprecht type, Maurey--Rosenthal coding and
rapidly increasing sequences based on $\ell_1$-averages---we add the
possibility of splitting an arbitrary vector into pieces of comparable norm,
while staying in one of these $\ell_\infty^n$'s.  This allows us to introduce
two additional classes of rapidly increasing sequences, and these in turn lead
to the stronger result about operators.

\subsection*{Acknowledgments}
Much of the research presented in this paper was carried out during
the second author's three visits to Athens in 2007 and 2008. He
offers his thanks to the National Technical University for the
support that made these visits possible and to all members of the
NTU Analysis group for providing an outstanding research
environment.  Both authors would especially like to thank  Ch.
Raikoftsalis for many stimulating discussions.

\section{Background}

\subsection{Notations}

We use standard notations: if $A$ is any set, $\ell_\infty(A)$ is the space of
all bounded (real-valued) functions on $A$, equipped with the supremum norm
$\|\cdot\|_\infty$ and $\ell_1(A) $ is the space of all absolutely summable
functions on $A$, equipped with the norm $\|x\|_1 = \sum_{a\in A}|x(a)|$.  The
{\em support} of a function $x$ is the set of all $a$ such that $x(a)\ne 0$;
$c_{00}(A)$ is the space of functions of finite support. We shall write
$\ell_p$ for the space $\ell_p(\nats)$, where $\nats$ is the set
$\{1,2,3,\dots\}$ of positive integers, and $\ell_p^n$ for
$\ell_p(\{1,2,\dots,n\})$.  Even when we are dealing with these sequence spaces
we shall use function notation $x(m)$, rather than subscript notation, for the
$m^{\text{th}}$ coordinate of the vector $x$.

When $x$ and $y$ are in $c_{00}(A)$ (and more generally) we shall write
$\langle y,x\rangle $ for $\sum_{a\in A} x(a)y(a)$. If we are thinking of $y$
as a functional acting on $x$ (rather than vice versa) we shall usually choose
a notation involving a star, denoting $y$ by $f^*$, or something of this kind.
In particular, $e_a$ and $e_a^*$ are two notations for the same unit vector in
$c_{00}(A)$ (given by $e_a(a')=\delta_{a,a'}$), to be employed depending on
whether we are thinking of it as a unit vector or as the evaluation functional
$x\mapsto \langle e^*_a, x\rangle = x(a)$. We apologize to those readers who
may find this kind of notation somewhat babyish.

We say that  (finitely or infinitely many)  vectors $z_1, z_2, \dots$ in
$c_{00}$ are {\em successive}, or that $(z_i)$ is a {\em block-sequence}, if
$\max \supp x_i<\min \supp x_{i+1}$ for all $i$. In a Banach space $X$ we say
that vectors $y_j$ are {\em successive linear combinations}, or that $(y_j)$ is
a {\em block sequence} of a basic sequence $(x_i)$ if  there exist
$0=q_1<q_2<\cdots$ such that, for all $j\ge 1$, $y_j$ is in the linear span
$[x_i :q_{j-1}<i\le q_j]$. If we may arrange that $y_j\in [x_i :q_{j-1}<i<
q_j]$ we say that $(y_j)$ is a {\em skipped block sequence.} More generally, if
$X$ has a Schauder decomposition $X=\bigoplus_{n\in \nats} F_n$ we say that
$(y_j)$ is a block sequence (resp. a skipped block sequence) with respect to
$(F_n)$ if there exist $0=q_0<q_1<\cdots$ such that $y_j$ is in
$\bigoplus_{q_{j-1}<n\le q_j}F_n$ (resp. $\bigoplus_{q_{j-1}<n< q_j}F_n$. A
{\em block subspace} is the closed subspace generated by a block sequence.

\subsection{Hereditary indecomposability}
A Banach space $X$ is {\em indecomposable} if there do not exist
infinite-dimensional closed subspaces $Y$ and $Z$ of $X$ with $X=Y\oplus Z$,
and is {\em hereditarily indecomposable} (HI) if every closed subspace is
indecomposable.  The following useful criterion, like so much else in this in
this area, goes back to the original paper of Gowers and Maurey \cite{GM}.

\begin{prop}\label{HICrit}
Let $X$ be a an infinite dimensional Banach space.  Then $X$ is HI if and only
if, for every pair $Y,Z$ of infinite-dimensional subspaces, and every
$\epsilon>0$, there exist $y\in Y$ and $z\in Z$ with $\|y+z\|>1$ and
$\|y-z\|<\epsilon$.  If $X$ has a finite-dimensional decomposition $(F_n)_{n\in
\nats}$ it is enough that the above should hold for block subspaces.
\end{prop}

We shall make use of the following well-known blocking lemma, the first part of
which can be found as Lemma~1 of \cite{M}. The proof of the second part is very
similar, and, as Maurey remarks, both can be traced back to R.C. James
\cite{J}.

\begin{lem}\label{ellOneAve}Let $n\ge 2$ be an integer, let
$\epsilon\in (0,1)$ be a real number and let $N$ be an integer that can be
written as $N=n^k$ for some $k\ge 1$. Let $(x_i)_{i=1}^N$ be a sequence of
vectors in the unit sphere of a Banach space $X$.

(i) If $\|\sum_{i=1}^N\pm x_i\|\ge (n-\epsilon)^k$ for all choices of signs
$\pm1$,  then there is a block sequence $y_1,y_2,\dots ,y_n\in [x_i:1\le i\le
N]$ which is $(1-\epsilon)^{-1}$-equivalent to the unit-vector basis of
$\ell_1^n$.

(ii) If $\|\sum_{i=1}^N\pm x_i\|\le (1+\epsilon)^k$ for all choices of signs
$\pm1$,  then there is a block sequence $y_1,y_2,\dots ,y_n\in [x_i:1\le i\le
N]$ which is $(1+\epsilon)$-equivalent to the unit-vector basis of
$\ell_\infty^n$.
\end{lem}

\subsection{$\mathscr L_{\infty}$-spaces}
A separable Banach space $X$ is an $\mathscr L_{\infty,\lambda}$-space if there
is an increasing sequence $(F_n)_{n\in \nats}$ of finite dimensional subspaces
of $X$ such that the union $\bigcup_{n\in \nats}F_n$ is dense in $X$ and, for
each $n$, $F_n$ is $\lambda$-isomorphic to $\ell_\infty^{\dim F_n}$. It is
known \cite{LS} that if a separable $\mathscr L_\infty$ space $X$ has no
subspace isomorphic to $\ell_1$, then the dual space $X^*$ is necessarily
isomorphic to $\ell_1$.  Of course, this implies that the dual of a separable,
hereditarily indecomposable $\mathscr L_\infty$-space is isomorphic to
$\ell_1$.

The Bourgain--Delbaen spaces $X_{a,b}$, which inspired the construction given
in this paper, were the first examples of $\mathscr L_{\infty}$ spaces not
containing $c_0$.

\subsection{Mixed Tsirelson spaces}\label{MT}

All existing HI constructions have, somewhere at the heart of them, a space of
Schlumprecht type; rather than working with the original space of \cite{S}, we
find it convenient to look at a different mixed Tsirelson space. We recall some
notation and terminology from \cite{AT}. Let $(l_j)_j$ be a sequence of
positive integers and let $(\theta_j)_j$ be a sequence of real numbers with
$0<\theta_n<1$. We define $W[(\mathscr A_{l_j},\theta_j)_j]$ to be the smallest
subset $W$ of $c_{00}$ with the following properties
\begin{enumerate}
\item $\pm e^*_k\in W$ for all $k\in \mathbb N$; \item whenever
$f^*_1,f^*_2,\dots,f^*_m\in W$ are successive vectors, $\theta_j\sum_{i\le m}
f^*_i \in W$, provided $m\le l_j$.
\end{enumerate}
We say that an element $f^*$ of $W$ is of Type 0 if $f^*=\pm e^*_k$ for some
$k$ and of Type I otherwise; an element of type I is said to have weight
$\theta_j$ if $f^*=\theta_j\sum_{i\le m} f^*_i $ for a suitable sequence
$(f_i)$ of successive elements of $W$.

The {\em  mixed Tsirelson space} $T[(\mathscr A_{l_j}, \theta_j)_j]$ is defined
to be the completion of $c_{00}$ with respect to the norm
$$
\|x\|= \sup \{ \langle f^*,x\rangle: f^*\in W[\mathscr A_{l_j}, \theta_j)_j].
$$
We may also characterize the norm of this space implicitly as being the
smallest function $x\mapsto \|x\|$ satisfying
$$
\|x\| = \max\bigg\{\|x\|_\infty, \sup \theta_j
\sum_{i=1}^{l_j}\|x\chi_{E_i}\|\bigg\},
$$
where the supremum is taken over all $j$ and all sequences of finite subsets
$E_1<E_2<\cdots<E_{l_j}$. Schlumprecht's original space is the result of taking
$l_j=j$ and $\theta_j = (\log_2(j+1))^{-1}$

In the rest of the paper we shall choose to work with two sequences of natural
numbers $(m_j)$ and $(n_j)$.  We require $m_j$ to grow quite fast, and $n_j$ to
grow even faster. The precise requirements are as follows.

\begin{assump}\label{mnAssump}
We assume that $(m_j,n_j)_{j\in \nats}$ satisfy the following:
\begin{enumerate}
\item $m_1\ge 4$; \item $m_{j+1} \ge m_j^2$; \item $n_{j+1} \ge
m_{j+1}^2(4n_j)^{2^{j+1}}.$
\end{enumerate}
\end{assump}

A straightforward way to achieve this is to assume that $(m_j,n_j)$ is some
subsequence of the sequence $(2^{2^j},2^{2^{j^2+1}})_{j\in \nats}$.  From now
on, whenever $m_j$ and $n_j$ appear, we shall assume we are dealing with
sequences satisfying \ref{mnAssump}.

 The following lemma can be found as
II.9 of \cite{AT}.  The proof is not affected by the small change we have made
in the definition of the sequences $(n_j)_j$ and $(m_j)_j$.

\begin{lem}\label{MTLem}
If $j_0\in \mathbb N$ and  $f\in W[(\mathscr A_{4n_j},m_j^{-1})_j]$ is an
element of weight $m_h$, then
$$
\bigg|\big\langle f^*,n_{j_0}^{-1}\sum_{l=1}^{n_{j_0}} e_l\big\rangle\bigg| \le
\begin{cases} 2m_h^{-1}m_{j_0}^{-1}\qquad\text{if } i<j_0\\
              m_h^{-1} \qquad\quad \text{if } i\ge j_0.
              \end{cases}
$$
In particular, the norm of $n_{j_0}^{-1}\sum_{l=1}^{n_{j_0}} e_l$ in
$T[(\mathscr A_{4n_j},m_j^{-1})_j]$ is exactly $m_{j_0}^{-1}$.

 If we restrict attention to $f\in
W[(\mathscr A_{4n_j},m_j^{-1})_{j\ne j_0}]$ then
$$
\bigg|\big\langle f^*,n_{j_0}^{-1}\sum_{l=1}^{n_{j_0}} e_l\big\rangle\bigg| \le
\begin{cases} 2m_h^{-1}m_{j_0}^{-2}\qquad\text{if } i<j_0\\
              m_h^{-1} \qquad \qquad\quad \text{if } i> j_0.
              \end{cases}
$$
In particular, the norm of $n_{j_0}^{-1}\sum_{l=1}^{n_{j_0}} e_l$ in
$T[(\mathscr A_{4n_j},m_j^{-1})_{j\ne j_0}]$ is at most $m_{j_0}^{-2}$.
\end{lem}

\section{The general Bourgain--Delbaen construction}
\label{BDSection}

In this section we shall present a generalization of the Bourgain--Delbaen
construction of separable $\mathscr L_\infty $-spaces.  Our approach is
slightly different from that of \cite{B} and \cite{BD}, but the mathematical
essentials are the same.  We choose to set things out in some detail partly
because we believe our approach yield new insights into the original BD
construction, and partly because the calculations presented here are a good
introduction to the notations and methods we use later. It is perhaps worth
emphasizing here that BD constructions are very different from the majority of
constructions that occur in Banach space theory. Normally we start with the
unit vectors in the space $c_{00}$ and complete with respect to some (possibly
exotic) norm.  The only norms that occur in a BD construction are the usual
norms of $\ell_\infty$ and $\ell_1$. What we construct here are exotic {\em
vectors} in $\ell_\infty$ whose closed linear span is the space we want.

The idea will be to introduce a particular kind of (conditional) basis for the
space $\ell_1$ and to study the subspace $X$ of $\ell_\infty$ spanned by the
biorthogonal elements.  Since $\ell_1$ is then in a natural way a subspace of
(and in some cases the whole of) $X^*$, we shall be thinking of elements of
$\ell_1$ as functionals and, in accordance with the convention explained
earlier, denote them $b^*$, $c^*$ and so on.  In our initial discussion we
shall consider the space $\ell_1(\nats)$ (which we shall later replace with
$\ell_1(\Gamma)$ with $\Gamma$ a certain countable set better adapted to our
needs).

\begin{defn} \label{TriangDef}We shall say that a basic sequence $(d_n^*)_{n\in
\nats}$ in $\ell_1(\nats)$ is a {\em triangular basis} if $\supp d^*_n\subseteq
\{1,2,\dots,n\}$, for all $n$.  We thus have
$$
d^*_n = \sum_{m=1}^n a_{n,m} e^*_m,
$$
where, by linear independence, we necessarily have $a_{n,n}\ne 0$. Notice that
the linear span $[d^*_1,d^*_2,\dots,d^*_n]$ is the same as
$[e^*_1,e^*_2,\dots,e^*_n]$, that is to say, the space $\ell_1(n)$, regarded as
a subspace of $\ell_1(\nats)$ in the usual way.  So, in particular, the basic
sequence $(d^*_n)$ is indeed a basis for the whole of $\ell_1$.  The
biorthogonal sequence in $\ell_\infty$ will be denoted $(d_n)$; it is a
weak$^*$ basis for $\ell_\infty$ and a basis for its closed linear span, which
will be our space $X$.
\end{defn}

\begin{prop}\label{TriangProp}
If $(d^*_n)$ is a triangular basis for $\ell_1(\nats)$, with basis constant
$M$, then the closed linear span $X=[d_n:n\in \nats]$ is a $\mathscr
L_{\infty,M}$-space.  If $(d^*_n)$ is boundedly complete, or equivalently
$(d_n)$ is shrinking, then $X^*$ is naturally isomorphic to $\ell_1(\nats)$
with $\|g^*\|_{X^*}\le \|g^*\|_1\le M\|g^*\|_{X^*}$.
\end{prop}
\begin{proof}
In accordance with our ``star'' notation, let us write $P^*_n$ for the basis
projection $\ell_1\to\ell_1$ associated with the basis $(d^*_n)$. Thus
$P^*_n(d^*_m)$ equals $d^*_m$ if $m\le n$ and $0$ otherwise; because $e^*_m\in
\ell_1(n)=[d^*_1,\dots,d^*_n]$, we also have $P^*_ne_m=e_m$ when $m\le n$. If
we modify $P^*_n$ by taking the codomain to be the image
$\text{im}\,P^*_n=\ell_1(n)$, rather than the whole of $\ell_1$, what we have
is a quotient operator, which we shall denote $q_n$, of norm at most $M$.  The
dual of this quotient operator is an isomorphic embedding
$i_n:\ell_\infty(n)\to \ell_\infty(\nats)$, also of norm at most $M$.  If $m\le
n$ and $u\in \ell_\infty(n)$ we have
$$
(i_nu)(m) = \langle e^*_m,i_nu\rangle= \langle q_ne^*_m,u\rangle = \langle
e^*_m,u\rangle = u(m).
$$
So $i_n$ is an {\em extension operator} $\ell_\infty^n\to \ell_\infty(\nats)$
and we have
$$
\|u\|_\infty \le \|i_nu\|_\infty \le M\|u\|_\infty
$$
for all $u\in \ell_\infty^n$.  In particular, the image of $i_n$, which is
exactly $[d_1,\dots, d_n]$ is $M$-isomorphic to $\ell_\infty^n$, which implies
that $X$ is a $\mathscr L_{\infty, M}$-space.

In the case where $(d^*_n)$ is a boundedly complete basis of $\ell_1$ then
$X^*$ may be identified with $\ell_1$ by standard result about bases.
Moreover, for $g^*\in \ell_1$, we have
$$
\|g^*\|_{X^*} = \sup \{\langle g^*,x\rangle: x\in X\text{ and }\|x\|_\infty\le
1\}\le \|g^*\|_1.
$$
On the other hand, if $g^*$ has finite support, say $\supp g^*\subseteq
\{1,2,\dots,n\}$, we can choose $u\in \ell_\infty^n$ with $\|u\|=1$ and
$\langle g^*,u\rangle =\|g^*\|_1$. The extension $x= i_n(u)$ is now in $X$ and
satisfies
$$
\|x\|\le M,\quad \langle g^*,x\rangle = \|g^*\|.
$$
\end{proof}

We shall say that $(d_n^*)$ is a {\em unit-triangular basis} of $\ell_1(\nats)$
if it is a triangular basis and the non-zero scalars $a_{n,n}$ are all equal to
1.  We can thus write
$$
d_n^* = e_n^* - c_n^*,
$$
where $c_1^*=0$ and $\supp c_n^*\subset \{1,2,\dots,n-1\}$ for $n\ge 2$.  The
clever part of the Bourgain--Delbaen construction is to find a method of
choosing the $c_n^*$ in such a way that $(d_n^*)$ is indeed a basic sequence.
The idea is to proceed recursively assuming that, for some $n\ge 1$, we already
have a unit-triangular basis $(d^*_m)_{m\le n}$ of $\ell_1^n$.  The value of
$P^*_rb^*$ is thus already determined when $1\le r\le n$ and $b^*\in \ell_1^n$.

\begin{defn} \label{BDDef}In the set-up described above, we shall say that an
element $c^*$ of $\ell_1(n)$ is a BD-{\em functional} (with respect to the
triangular basis $(d^*_m)_{m=1}^n$) if there there exist real numbers
$\alpha\in (0,1]$ and $\beta\in [0,\frac12)$ such that we can express $c^*$ in
one of the following forms:
\begin{enumerate}
\setcounter{enumi}{-1}
 \item $\alpha e^*_j$ with $1\le j\le n$,
 \item $\beta (I-P^*_k) b^*$ with $0\le k<n$ and $b^*\in
\ball \ell_1(k+1,\dots,n)$,
 \item $\alpha e^*_j+\beta (I-P^*_k)b^*$ with $1\le j\le k<n$ and $b^*\in
\ball \ell_1(k+1,\dots,n)$.
 \end{enumerate}
 The non-negative constant $\beta$ will
be called the {\em weight} of the functional $c^*$ (``weight 0'' in case (0)).
Note that (0) and (1) are ``almost''  special cases of (2), with $\beta$ (resp.
$\alpha$) equal to 0.  In the construction presented in this paper, we do not
use functionals of type (0) and the constant $\alpha$ in case (2) is always
equal to 1. However, it may be worth stating the following theorem in full
generality.
\end{defn}

\begin{thm}[\cite{B}, \cite{BD}] \label{BDThm}
Let $\theta$ be a real number with $0<\theta <\frac12$ and let $d_n^* =
e^*_n-c^*_n$ in $\ell_1$ be such that, for each $n$, $c_{n+1}^*\in \ell_1^n$ is
a BD-functional of weight at most $\theta$ with respect to $(d^*_m)_{m=1}^n$.
Then $(d^*_n)_{n\in \nats}$ is a triangular basis of $\ell_1$, with basis
constant at most $M= 1/(1-2\theta)$.  The subspace $X=[d_n:n\in \nats]$ of
$\ell_\infty$ is thus a $\mathscr L_{\infty, M}$-space.
\end{thm}
\begin{proof} Despite the disguise, this is essentially the same
argument as in the original papers of Bourgain and Delbaen. What  we need to
show is that $P^*_m$ is a bounded operator, with $\|P^*_m\|\le M$ for all $m$.
Because we are working on the space $\ell_1$ it is enough to show that
$\|P^*_me^*_n\|\le M$ for every $m$ and $n$.

First, if $n\le m$, $P^*_me^*_n=e^*_n$, so there is nothing to prove. Now  let
us assume that $\|P_ke^*_j\|\le M$ for all $k\le m$ and all $j\le n$; we then
consider $P^*_me^*_{n+1}$. We use the fact that
$$
e^*_{n+1}= d^*_{n+1}+ c_{n+1}^*, $$
 with $c_{n+1}^*\in \ell_1^n$ a
BD-functional.  We shall consider a functional of type (2), which presents the
most difficulty. We thus have
$$
c_{n+1}^* = \alpha e_j^* + \beta (I-P^*_k)b^*,
$$
where $1\le j\le k< n$ and $\alpha, \beta, b^*$ are as in
Definition~\ref{BDDef}, and $\beta\le \theta$ by our hypothesis.
Now, because $n+1>m$ we have $P^*_md^*_{n+1}=0$ so
$$
P^*_me^*_{n+1} = \alpha P^*_me^*_j + \beta (P^*_m - P^*_{m\wedge k})b^*.
$$

If $k\ge m$ the second term vanishes so that
$$
\|P^*_me^*_n\|= \alpha\|P^*_me^*_j\| \le \|P^*_me^*_j\|,
$$
which is at most $M$ by our inductive hypothesis.

If, on the other hand, $k<m$, we certainly have $j<m$ so that
$P_m^*e^*_j=e^*_j$, leading to the estimate
$$
\|P^*_me^*_{n+1}\|\le  \alpha\|e^*_j\| + \beta\|P^*_mb^*\| + \beta\|P^*_kb^*\|.
$$
Now $b^*$ is a convex combination of functionals $\pm e^*_l$ with $l\le n$, and
our inductive hypothesis is applicable to all of these.  We thus obtain
$$
\|P^*_me^*_{n+1}\|\le  \alpha + M\beta\le 1+2M\beta= M,
$$
by the definition of $M=1/(1-2\theta)$ and the assumption that $0\le\beta\le
\theta$.
\end{proof}

The $\mathscr L_\infty$ spaces of Bourgain and Delbaen, and those we construct
in the present paper are of the above type.  However, the ``cuts'' $k$ that
occur in the definition of $BD$-functionals are restricted to lie in a certain
subset of $\nats$, thus naturally dividing the coordinate set $\nats$ into
successive intervals. As in \cite{H}, it will be convenient to replace the set
$\nats$ with a different countable set $\Gamma$ having a structure that
reflects this decomposition. This will also enable us later to use a notation
in which an element $\gamma\in \Gamma$ automatically codes the BD-functional
associated with it.

\begin{thm}\label{BDThmBis}
Let $(\Delta_q)_{q\in \nats}$ be a sequence of non-empty finite sets, with
$\#\Delta_1=1$; write $\Gamma_q=\bigcup_{1\le p\le q}\Delta_p$,
$\Gamma=\bigcup_{ p\in \nats}\Delta_p$. Assume that there exists
$\theta<\frac12$ and a mapping $\tau$ defined on $\Gamma\setminus \Delta_1$,
assigning to each $\gamma\in \Delta_{q+1}$ a tuple of one of the forms:
\begin{enumerate}
\setcounter{enumi}{-1}
 \item $(\alpha, \xi)$ with $0<\alpha\le 1$ and $\xi\in \Gamma_q$;
 \item $(p,\beta, b^*)$ with $0\le p< q$, $0<\beta\le \theta$ and $b^*\in \ball
 \ell_1\left(\Gamma_q\setminus \Gamma_p\right)$;
 \item $(\alpha,\xi,p,\beta,b^*)$ with $0<\alpha\le 1$, $1\le p <q$, $\xi\in
 \Gamma_p$, $0<\beta\le \theta$  and $b^*\in \ball
 \ell_1\left  (\Gamma_q\setminus \Gamma_p\right)$.
 \end{enumerate}
Then there exist $d_\gamma^* = e^*_\gamma-c^*_\gamma\in \ell_1(\Gamma)$ and
projections $P^*_{(0,q]}$ on $\ell_1(\Gamma)$ uniquely determined by the
following properties:
\begin{enumerate}
  \item
  $\displaystyle P^*_{(0,q]}d^*_\gamma = \begin{cases}
 d^*_\gamma\qquad\qquad\qquad\qquad\qquad\text{if }\gamma\in \Gamma_q\\
 0\ \qquad\qquad\qquad\qquad\qquad\text{if }\gamma\in \Gamma\setminus
 \Gamma_q
 \end{cases}$
 \item $\displaystyle c^*_\gamma =\qquad \begin{cases}
  0  \qquad\qquad\qquad\qquad\qquad\text{if }\gamma\in \Delta_1\\
  \alpha e^*_\xi \qquad \qquad\qquad\text{if } \tau(\gamma) = (\alpha,
  \xi)\\
  \beta (I-P^*_{(0,p]}) b^*\qquad\ \text{
  if }\tau(\gamma)=(p, \beta, b^*)\\
        \alpha e^*_\xi + \beta(I-P^*_{(0,p]}) b^*\quad\text{
 if }\tau(\gamma)=(\alpha, \xi, \beta, b^*)\text{ with } \xi\in \Delta_p.
 \end{cases}$
 \end{enumerate}
 The family $(d^*_\gamma)_{\gamma\in \Gamma)}$ is a basis for
 $\ell_1(\Gamma)$ with basis constant at most $M=
 (1-2\theta)^{-1}$.  The norm of each projection $P^*_{(0,q]}$ is at
 most $M$.  The biorthogonal elements $d_\gamma$ generate a
 $\mathscr L_{\infty ,M}$-subspace $X(\Gamma,\tau)$ of $\ell_\infty(\Gamma)$.
 For each $q$ and each $u\in \ell_\infty(\Gamma_q)$, there is a
 unique $i_q(u)\in [d_\gamma:\gamma\in \Gamma_q]$ whose restriction
 to $\Gamma_q$ is $u$; the extension operator
 $i_q:\ell_\infty(\Gamma_q)\to X(\Gamma,\tau)$ has norm at most $M$.
  The subspaces $M_n=[d_\gamma:\gamma\in \Delta_q]
  =i_q[\ell_\infty(\Delta_q)]$ form a
finite-dimensional decomposition (FDD) for $X$; if this FDD is shrinking then
$X^*$ is naturally isomorphic to $\ell_1(\Gamma)$.
\end{thm}
\begin{proof} We shall show that, with a suitable identification of $\Gamma$ with
$\nats$, this theorem is just a special case of Theorem~\ref{BDThm}.
Let $k_p=\#\Gamma_p$ and let $n\mapsto \gamma(n):\nats\to \Gamma$ be
a bijection with the property that $\Delta_1=\{\gamma(1)\}$, while,
for each $q\ge2$, $\Delta_q = \{\gamma(n):k_{q-1}<n\le k_q\}$. There
is a natural isometry: $J:\ell_1(\nats)\to \ell_1(\Gamma)$
satisfying $J(e^*_n)=e^*_{\gamma(n)}$.  It is straightforward to
check that if $d^*_n=J^{-1}(d^*_{\gamma(n)})=e^*_n-c^*_n$, then the
hypotheses of Theorem~\ref{BDThm} are satisfied. (The cuts $k$ that
occur in the BD-functionals $c^*_n$ are all of the form $k=k_p$.)
All the assertions in the present theorem are now immediate
consequences. The projections $P^*_{(0,q]}$ whose existence is
claimed here are given by $P^*_{(0,q]}= JP^*_{k_q}J^{-1}$, where
$P^*_n$ is the basis projection of Theorem~\ref{BDThm}. When ordered
as $(d_{\gamma(n)})_{n\in \nats}$ the vectors $d_\gamma$ form a
basis of their closed linear span, which is a $\mathscr L_{\infty,
M}$-space. The extension operator that (by abuse of notation) we
here denote by $i_q$ is just $Ji_{k_q}J^{-1}$. The assertions about
the subspaces  $M_q=[d_{\gamma(n)}:k_{q-1}<n\le k_q]$ follow from
the fact that $(d_{\gamma(n)})$ is a basis.
\end{proof}

We now make a few observations  about the space $X=\mathfrak(\Gamma,\tau)$ and
the functions $d_\gamma$, taking the opportunity to introduce notation that
will be used in the rest of the paper. We have seen that for each $\gamma\in
\Delta_{n+1}$ the functional $d^*_\gamma$ has support contained in
$\Gamma_{n}\cup\{\gamma\}$. Using biorthogonality, we see that $d_\gamma$ is
supported by $\{\gamma\} \cup \Gamma\setminus \Gamma_{n+1}$.  It should be
noted that we should not expect the support of $d_\gamma$ to be finite; in
fact, in all interesting cases, we  have $X\cap c_0(\Gamma)=\{0\}$.

As noted above the subspaces $M_n = [d_\gamma:\gamma\in \Delta_n]$
form a finite-dimensional decomposition for $X$. For each interval
$I\subseteq \nats$ we define the projection $P_I:X\to \bigoplus
_{n\in I}M_n$ in the natural way; this is consistent with our use of
$P^*_{(0,n]}$ in Theorem~\ref{BDThmBis}. Most of our arguments will
involve sequences of vectors that are block sequences with respect
to this FDD.  Since we are using the word ``support'' to refer to
the set of $\gamma$ where a given function is non-zero, we need
other terminology for the set of $n$ such that $x$ has a non-zero
component in $M_n$. We define the {\em range}  of $x$, denoted $\ran
x$, to be the smallest interval $I\subseteq \nats$ such that $x\in
\bigoplus_{n\in I}M_n$.  It is worth noting that if $\ran x=(p,q]$
then we can write $x=i_q(u)$ where $u=x\restriction \Gamma_q\in
\ell_\infty(\Gamma_q)$ satisfies $\Gamma_p\cap\supp u=\emptyset.$

\section{Construction of $\BmT$
and $\XK$}

We now set about constructing specific BD spaces which will be
modelled on mixed Tsirelson spaces, in rather the same way that the
original spaces of Bourgain and Delbaen have been found to be
modelled on $\ell_p$. We shall adopt a notation in which elements
$\gamma$ of $\Delta_{n+1}$ automatically code the corresponding
BD-functionals. This will allow us to write $X(\Gamma)$ rather than
$X(\Gamma,\tau)$ for the resulting $\mathscr L_\infty$-space. To be
more precise, an element $\gamma$ of $\Delta_{n+1}$ will be a tuple
of one of the forms:
\begin{enumerate}
\item $\gamma= (n+1,\beta,b^*)$,\quad in which case $\tau(\gamma)=(0,\beta,b^*)$;
\item $\gamma = (n+1,\xi,\beta,b^*)$ in which case $\tau(\gamma)=(1,\xi,\rank \xi, \beta,b^*)$.
\end{enumerate}
In each case, the first co-ordinate of $\gamma$ tells us what the {\em rank} of
$\gamma$ is, that is to say to which set $\Delta_{n+1}$ it belongs, while the
remaining co-ordinates specify the corresponding BD-functional.

It will be observed that BD-functionals of Type 0 do not arise in this
construction and that the $p$ in the definition of a Type 1 functional is
always 0. In the definition  of a Type 2 functional that the scalar $\alpha$
that occurs
 is always 1 and $p$ equals $\rank \xi$. We shall make
the further restriction the weight $\beta$ must be of the form $m_j^{-1}$,
where the sequences $(m_j)$ and $(n_j)$ satisfy Assumption~\ref{mnAssump}.  We
shall say that the element $\gamma$ has {\em weight} $m_j^{-1}$ (sometimes
dropping the $^{-1}$ and referring to ``weight $m_j$''). In the case of a Type
2 element $\gamma=(n+1, \xi, m^{-1}_j,b^*)$ we shall insist that $\xi$ be of
the same weight $m_j^{-1}$ as $\gamma$.

To ensure that our sets $\Delta_{n+1}$ are finite we shall admit into
$\Delta_{n+1}$ only elements of weight $m_j$ with $j\le n+1$. A further
restriction involves a recursively defined function which we call ``age''. For
a Type 1 element $\gamma=(n+1, \beta, b^*)$ we define $\age\gamma=1$. For a
Type 2 element $\gamma=(n+1, \xi, m^{-1}_j,b^*)$, we define $\age \gamma= 1 +
\age \xi$, and further restrict the elements of $\Delta_{n+1}$ by insisting
that the age of an element of weight $m_j$ may not exceed $n_j$. Finally, we
shall restrict the functionals $b^*$ that occur in an element of $\Delta_{n+1}$
by requiring them to lie in some finite subset $B_n$ of $\ell_1(\Gamma_n)$. It
is convenient to fix an increasing sequence of natural numbers $(N_n)$ and take
$B_{n,p}$ to be the set of all linear combinations $b^*=\sum_{\eta\in
\Gamma_n\setminus \Gamma_p}a_\eta e^*_\eta$, where $\sum_\eta|a_\eta|\le 1$ and
each $a_\eta$ is a rational number with denominator dividing $N_n!$. We may
suppose the $N_n$ are chosen in such a way that $B_{n,p}$ is a $2^{-n}$-net in
the unit ball of $\ell_1(\Gamma_n\setminus \Gamma_p)$. The above restrictions
may be summarized as follows.

\begin{assump}\label{DeltaUpperAssump}
\begin{align*}
\Delta_{n+1} &\subseteq \bigcup_{j=1}^{n} \left\{(n+1,m_j^{-1},b^*): b^*\in B_{n,0}\right\}\\
&\cup  \bigcup_{0<p<n}\bigcup_{j=1}^{p}\left\{(n+1,\xi,m_j^{-1},b^*): \xi\in
\Delta_p, \weight \xi = m_j^{-1},\ \age\xi<n_j,\ b^*\in B_{n,p}\right\}
\end{align*}
\end{assump}

We shall also assume that $\Delta_{n+1}$ contains a rich supply of elements of
``even weight'', more exactly of weight $m_j$ with $j$ even.

\begin{assump}\label{DeltaLowerAssump}
\begin{align*}
\Delta_{n+1} &\supseteq \bigcup_{j=1}^{\lfloor( n+1)/2\rfloor} \left\{(n+1,m_{2j}^{-1},b^*): b^*\in B_{n,0}\right\}\\
&\cup\bigcup_{1\le p<n} \bigcup_{j=1}^{\lfloor p/2\rfloor}
\left\{(n+1,\xi,m_{2j}^{-1},b^*): \xi\in \Delta_p, \weight \xi =
m_{2j}^{-1},\ \age\xi<n_{2j},\ b^*\in B_{n,p}\right\}
\end{align*}
\end{assump}

For our main HI construction, there are additional restrictions on the elements
with ``odd weight'' $m_{2j-1}$.  However, there is some interest already in the
space we obtain without making such restrictions. We denote this space $\BmT$;
it is an isomorphic predual of $\ell_1$ that is unconditionally saturated but
contains no copy of $c_0$ or $\ell_p$.  An analogous space $\mathfrak
B_{\text{T}}$, modelled on the standard Tsirelson space, rather than a mixed
Tsirelson space, was constructed a few years ago by the second-named author.

\begin{defn}
We define $\BmT=\mathfrak B_{\text{mT}}[(m_j,n_j)_{j\in \nats}]$ to be the
space $X(\Gamma)$ where $\Gamma=\Gamma^{\text{max}}$ is defined by the
recursion $\Delta_1=\{1\}$,
\begin{align*}
\Delta_{n+1} &=\bigcup_{j=1}^{n+1} \left\{(n+1,m_j^{-1},b^*): b^*\in B_{n,0}\right\}\\
&\cup \bigcup_{j=1}^{n-1} \bigcup_{j\le p<n}\left\{(n+1,\xi,m_j^{-1},b^*):
\xi\in \Delta_p, \weight \xi = m_j^{-1},\ \age\xi<n_j,\ b^*\in B_{n,p}\right\}
\end{align*}
\end{defn}

The extra constraints that we place on  ``odd-weight'' elements in order to
obtain hereditary indecomposability will involve a coding function that will
produce the analogues of the ``special functionals'' that occur in \cite{GM}
and other HI constructions. In our case, all we need is an injective function
$\sigma:\Gamma\to \nats$ satisfying $4\sigma (\gamma)>\rank \gamma$ for all
$\gamma$. This may easily be included in our recursive construction of
$\Gamma$.  We then insist that a Type 1 element of odd weight must have the
form
$$
(n+1, m_{2j-1}^{-1}, e^*_{\eta})
$$
with $\weight \eta = m_{4i-2}>n_{2j-1}^2$, while a Type 2 element must be
$$
(n+1,\xi,m_{2j-1}^{-1}, e^*_\eta)
$$
with $\weight \eta = m_{4\sigma(\xi)}$.

\begin{defn}\label{XKDef}
We define $\mathfrak X_{\text{K}}[(m_j,n_j)_{j\in \nats}]$ to  be the space
$X(\Gamma)$ where $\Gamma=\Gamma^{\text{K}}$ is defined by the recursion
$\Delta_1=\{1\}$,
\begin{align*}
\Delta_{n+1} &= \bigcup_{j=1}^{\lfloor(n+1)/2\rfloor}
\left\{(n+1,m_{2j}^{-1},b^*): b^*\in B_{n,0}\right\}\\ &\cup
\bigcup_{p=1}^{n}\bigcup_{j=1}^{\lfloor
p/2\rfloor}\left\{(n+1,\xi,m_{2j}^{-1},b^*): \xi\in \Delta_p, \weight \xi =
m_{2j}^{-1},\ \age\xi<n_{2j},\
b^*\in B_{n,p}\right\}\\
&\cup\bigcup_{j=1}^{\lfloor(n+2)/2\rfloor}
\left\{(n+1,m_{2j-1}^{-1},e^*_\eta):
\eta\in \Gamma_n \text{ and } \weight \eta= m_{4i-2}>n_{2j-1}^2\right\}\\
&\cup  \bigcup_{p=1}^{n}\bigcup_{j=1}^{\lfloor(
p+1)/2\rfloor}\left\{(n+1,\xi,m_{2j-1}^{-1},e^*_\eta): \xi\in \Delta_p, \weight
\xi = m_{2j-1}^{-1}\right.,\\  &\left.
\qquad\qquad\qquad\qquad\qquad\qquad\age\xi<n_{2j-1},\ \eta\in
\Gamma_n\setminus\Gamma_p,\ \weight\eta= m_{4\sigma(\xi)}\right\}.
\end{align*}
\end{defn}

With the definition readily at hand, this is a convenient moment to record an
important ``tree-like'' property of odd-weight elements of $\Gamma^{\mathrm
K}$, even though we shall not be exploiting these special elements until later
on.

\begin{lem}\label{Treelike}
Let $\gamma, \gamma'$ be two elements of $\Gamma^{\mathrm K}$ both of weight
$m_{2j-1}$ and of ages $a\ge a'$, respectively.  Let
$(p_i,e^*_{\eta_i},\xi_i)_{1\le i\le a}$, resp.
$(p'_i,e^*_{\eta'_i},\xi'_i)_{1\le i\le a'}$, be the analysis of $\gamma$,
resp. $\gamma'$. There exists $l$ with $1\le l\le a'$ such that $\xi'_i=\xi_i$
when $i<l$, while $\weight \eta_j \ne \weight \eta'_i$ for all $j$ when $l<i\le
a'$.
\end{lem}
\begin{proof}
If $\weight\eta'_i\ne\weight\eta_j$ for all $i\ge 2$ and all $j$ there is
nothing to prove (we may take $l=1$).  Otherwise, let $2\le l\le a$ be maximal
subject to the existence of $j$ such that $\weight\eta_j=\weight\eta'_l$.  Now
this weight is exactly $m_{4\sigma(\xi'_{l-1})}$, which means that $j$ cannot
be 1 (because the weight of $\eta_1$ has the form $m_{4k-2}$).  Thus
$\sigma(\xi'_{l-1})=\sigma(\xi_{j-1})$, which implies that
$\xi'_{l-1}=\xi_{j-1}$.  Since $l-1=\age\xi'_{l-1}$ and $j-1=\age\xi_{j-1}$, we
deduce that $j=l$.  Moreover, since the elements $\xi_{i}$ with $i<l-1$ are
determined by $\xi_{l-1}$, we have $\xi_i=\xi'_i$ for $i<l$.
\end{proof}

Although the structure of the space $X(\Gamma)$ is most easily understood in
terms of the basis $(d_\gamma)$ and the biorthogonal functionals $d_\gamma^*$,
it is with the evaluation functionals $e^*_\gamma$ that we have to deal in
order to estimate norms.  The recursive definition of the functionals
$d_\gamma^*$ can be unpicked to yield the following proposition.

\begin{prop} \label{EvalAnal}Assume that the set $\Gamma$ satisfies
Assumption~\ref{DeltaUpperAssump}. Let $n$ be a positive integer and let
$\gamma$ be an element of $\Delta_{n+1}$ of weight $m_j$ and age $a\le n_j$.
Then there exist natural numbers $0=p_0<p_1<\cdots<p_a=n+1$, elements
$\xi_1,\dots,\xi_a=\gamma$ of weight $m_j$ with $\xi_r\in \Delta_{p_r}$ and
functionals $b^*_r\in \ball\ell_1\left(\Gamma_{p_r-1}\setminus
\Gamma_{p_{r-1}}\right)$ such that
 \begin{align*} e^*_\gamma &= \sum_{r=1}^a d^*_{\xi_r} +
m_j^{-1}\sum_{r=1}^a P^*_{(p_{r-1},\infty)}b^*_r\\
&=\sum_{r=1}^a d^*_{\xi_r} + m_j^{-1}\sum_{r=1}^a P^*_{(p_{r-1},p_r)}b^*_r.
\end{align*}
\end{prop}
\begin{proof}
Given the assumption \ref{DeltaUpperAssump}, this is an easy induction on the
age $a$ of $\gamma$. If $a=1$ then $\gamma$ has the form $(n+1, m_j^{-1}, b^*)$
and
$$
e^*_\gamma = d^*_\gamma + c^*_\gamma,
$$
where $c^*_\gamma$ is the Type 1 BD-functional
$$
c^*_\gamma = m_j^{-1}P^*_{(0,\infty)}b^*,
$$
with $b^*\in B(n,0)\subset \ball \ell_1\left(\Gamma_n\right)$. Since $b^*$ is
in the image of the projection $P^*_{(0,n]}$ we have $P^*_{(0,n]}b^*=b^*$ and
so
$$
e^*_\gamma = d^*_{\xi_1} + m_j^{-1}P^*_{(p_0,\infty)}b^*_1=d^*_{\xi_1} +
m_j^{-1}P^*_{(p_0,p_1)}b^*_1,
$$
with $p_0=0$, $p_1=n+1$, $b^*_1=b^*$ and $\xi_1=\gamma$.

If $a>1$ then $\gamma$ has the form $(n+1,\xi,m_j^{-1},b^*)$ and $c^*_\gamma$
is the Type 2 BD-functional
$$
c^*_\gamma = e^*_\xi + m_j^{-1} P^*_{(p,\infty)}b^*.
$$
If we apply our inductive hypothesis to the element $\xi$ of weight $m_j$, rank
$p$ and age $a-1$, we obtain the desired expression for $e^*_\gamma$.
\end{proof}

We shall refer to the identity presented in the above proposition as the {\em
evaluation analysis} of $\gamma$ and shall use it repeatedly in norm
estimations. The form of the second term in the evaluation analysis, involving
a sum weighted by $m_j^{-1}$, indicates that there is going to be a connection
with mixed Tsirelson spaces; the first term, involving functionals $d^*_\xi$,
with no weight, can cause inconvenience in some of our calculations, but is an
inevitable feature of the BD construction. The data $(p_r,b^*_r,\xi_r)_{1\le
r\le a}$ will be called the {\em analysis} of $\gamma$. We note that if $1\le
s\le a$ the analysis of $\xi_s$ is just $(p_r,b^*_r,\xi_r)_{1\le r\le s}$.

In the remainder of this section, and in the next, we shall be dealing with a
space $X=X(\Gamma)$ and shall be making the assumptions \ref{DeltaUpperAssump}
and \ref{DeltaLowerAssump}.  Our results thus apply both to $\BmT$ and $\XK$.
We note that, since the weights $m_j^{-1}$ are all at most $\frac14$, the
constant $M$ in Theorem~\ref{BDThmBis} may be taken to be 2.  This leads to the
following norm estimates for the extension operators $i_n$ and for the
projections $P_I$ associated with the FDD $(M_n)$:
$$
\|i_n\|=\|P_{(0,n]}\|\le 2,\quad \|P_{(n,\infty)}\|\le 3, \quad
\|P_{(m,n]}\|\le 4, \quad \|d^*_\xi\|=\|P^*_{[\rank\xi,\infty)}e^*_\xi\|\le 3.
$$
The assumption \ref{DeltaLowerAssump} enables to write down a kind of converse
to Proposition~\ref{EvalAnal} which will lead to our first norm estimate.

\begin{prop} \label{ExistGammaEven}Let $j,a$ be positive integers with $a\le n_{2j}$, let
 $0=p_0<p_1<p_2<\cdots<p_a$ be natural numbers with $p_1\ge 2j$ and
let $b^*_r$ be functionals in $B(p_r-1,p_{r-1})$ for $1\le r\le a$. Then there
are elements $\xi_r\in \Gamma_{p_r}$ such that the analysis of $\gamma=\xi_a$
is $(p_r,b^*_r,\xi_r)_{1\le r\le a}$.
\end{prop}
\begin{proof}
This is another easy induction on $a$.  For $a=1$, the assumption
\ref{DeltaLowerAssump} and the hypothesis that $p_1\ge 2j$ guarantee that the
tuple $\xi_1=(p_1,m_{2j}^{-1},b^*_1)$ is in $\Gamma_{p_1}$. We continue
recursively, setting $\xi_{r+1} = (p_{r+1}, \xi_r, m_{2j}^{-1}, b^*_{r+1})$.
\end{proof}

\begin{prop}\label{LowerEst} Let $(x_r)_{r=1}^a$ be a skipped block
sequence (with respect to the FDD $(M_n)$) in $X$.  If $j$ is a positive
integer such that $a\le n_{2j}$ and $2j<\min\ran x_2$, then there exists an
element $\gamma$ of weight $m_{2j}$ satisfying
\begin{align*}
\sum_{r=1}^a x_r(\gamma)&\ge \half m_{2j}^{-1} \sum_{r=1}^a
\|x_r\|.\\
\intertext{Hence} \|\sum_{r=1}^a x_r\|&\ge \half m_{2j}^{-1} \sum_{r=1}^a
\|x_r\|.
\end{align*}
\end{prop}
\begin{proof}
Let $p_0=0$,  and choose  $p_1,p_2,\dots,p_a$ such that $\ran
x_r\subseteq (p_{r-1},p_r)$. Thus $x_r = i_{p_r-1}(u_r)$ where the
element $u_r=x_r\restriction \Gamma_{p_r-1}$ has support disjoint
from $\Gamma_{p_{r-1}}$.  Since $\|i_n\|\le 2$ for all $n$ we have
$\|u_r\|\ge\frac12 \|x_r\|$ and so there exist $\eta_r\in
\Gamma_{p_r-1}\setminus \Gamma_{p_{r-1}}$ with
$$
|u_r(\eta_r)|\ge \frac12 \|x_r\|.
$$
The functional $b^*_r=\pm e^*_{\eta_r}$ is certainly in $B_{p_{r}-1,p_{r-1}}$
and  with a suitable choice of sign we may arrange that
$$
\langle b^*_r, x_r\rangle = |u_r(\eta_r)|\ge  \frac12 \|x_r\|.
$$

By Proposition~\ref{ExistGammaEven} there is an element $\gamma$ of
$\Delta_{p_a}$ whose analysis is $(p_r,b^*_r,\xi_r)_{1\le r\le a}$. We shall
use the evaluation analysis to calculate
$$
\sum_{s=1}^a x_s(\gamma)= \langle e^*_\gamma, \sum_{s=1}^a x_s\rangle.
$$
For any $r$ and $s$, $x_s\in [d_\xi:p_{s-1}<\rank \xi<p_s]$, while $\rank
\xi_r=p_r$, whence
$$
\langle d^*_{\xi_r},x_s\rangle = 0\text{ \quad for all } r,s,
$$
while
 $$
\langle   P_{(p_{r-1},p_r)}^*b^*_r, x_s\rangle = \langle b^*_r,
P_{(p_{r-1},p_r)}x_s\rangle = 0,
  $$
for all $r\ne s$.  In the case $r=s$ we have
$$
\langle   P_{(p_{r-1},p_r)}^*b^*_r, x_r\rangle = \langle b^*_r,
P_{(p_{r-1},p_r)}x_r\rangle = \langle b^*_r,x_r\rangle.
  $$

The evaluation analysis thus simplifies to yield
$$
\sum_{r=1}^a x_r(\gamma) = m_{2j}^{-1} \sum_{r=1}^a \langle
b^*_r,x^*r\rangle\ge \half m_{2j}^{-1}\sum_{r=1}^a \|x_r\|.
$$
\end{proof}

The lower estimate we have just obtained indicates that there is a close
connection between our space $X$ and mixed Tsirelson spaces of the kind
considered in Subsection~\ref{MT}. With a bit more work one can show that a
normalized skipped-block sequence in $X$ dominates the unit vector basis of
$T[(\mathscr A_{n_{2j}}, m_{2j}^{-1})_{j\in \nats}]$.  We shall not need this
more precise result in the present work.

\section{Rapidly increasing sequences}

We continue  to work with the space $X=X(\Gamma)$, where $\Gamma$ satisfies the
assumptions \ref{DeltaUpperAssump} and \ref{DeltaLowerAssump}. We saw in the
last section that skipped block sequences admit useful Mixed Tsirelson lower
estimates.  We now pass to a class of block sequences that admit upper
estimates of a similar kind. The following definition is a variant of something
that is familiar from other HI constructions.

\begin{Defn}\label{RISDef}
Let $I$ be an interval in $\nats$ and let $(x_k)_{k\in I}$ be a block sequence
(with respect to the FDD $(M_n)$). We say that $(x_k)$ is a {\em rapidly
increasing sequence}, or RIS, if there exists a constant $C$ such that the
following hold:
\begin{enumerate}
\item
$\|x_k\|\le C$ for all $k\in \nats$,
\end{enumerate} and there is an increasing sequence $(j_k)$ such
that, for all $k$,
\begin{enumerate}\setcounter{enumi}{1}
\item $j_{k+1} > \max\, \text{ran }\,x_k$
\item $|x_k(\gamma)| \le Cm_i^{-1}$ whenever $\weight\gamma = m_i$ and
$i<j_k$
\end{enumerate}
If we need to be specific about the constant, we shall refer to a sequence
satisfying the above conditions as a $C$-RIS.
\end{Defn}

\begin{lem}\label{RISEst} Let $(x_k)$ be a $C$-RIS and let $(j_k)$ be an
increasing sequence of natural numbers as in the definition. If $\gamma\in
\Gamma$ and  $\weight\gamma= m_i$ then, for any natural number $s$
$$
|\langle e_\gamma^*, P_{(s,\infty)}x_k\rangle| \le
\begin{cases}5Cm_i^{-1}\quad\text{ if }i<j_k \\
 3Cm_i^{-1} \quad\text{ if }i\ge  j_{k+1}.\end{cases}
$$
\end{lem}
\begin{proof}
We first consider the case where $i\ge j_{k+1}$, noting that this implies that
$i>\max\, \text{ran }\,x_k$ by RIS condition (2). As in
Proposition~\ref{EvalAnal}, we may write down the evaluation analysis of
$\gamma$ as
$$
e_\gamma^*= \sum_r d^*_{\xi_r} + m_i^{-1} \sum_r b^*_r\circ
P_{(p_{r-1},\infty)},
$$
where $ 0=p_0<p_1<q_1< p_2<\cdots, $ and $b^*_r$ is a norm-1 element of
$\ell_1(\Gamma)$, supported by $\Gamma_{p_r-1}\setminus \Gamma_{p_{r-1}}$,
while $\xi_r$ is of rank $p_r$ and weight $m_i$. Since $\Delta_{q}$ contains no
elements of weight $m_i$ unless $q\ge i$, it must be that $p_1\ge i$.  Thus
$p_1> \max\, \text{ran }\,x_k$, from which it follows that
$P_{(p_r,\infty)}\circ P_{(s,\infty)}x_k=P_{(s\vee p_r,\infty)}x_k=0$ for all
$r\ge 1$. For the same reason, we also have
 $$\langle d^*_{\xi_r},
P_{(s,\infty)}x_k\rangle =\langle e^*_{\xi_r}, P_{(s\vee
q_r,\infty)}P_{[p_r,\infty )}x_k=0
$$
for all $r$. We are left with
$$
|\langle e_\gamma^*, P_{(s,\infty)}x_k\rangle| = m_i^{-1}|\langle
b^*_1, P_{(s,\infty)}x\rangle| \le m_i^{-1}\|P_{(s,\infty)}\|\,\|x\|
\le 3C m_i^{-1}
$$

In the case where $i<j_k$, we again use the evaluation analysis, but need to be
more careful about the value of $s$. Since we shall need this argument again,
we state it as a separate lemma. Clearly the second part of the present lemma
is an immediate consequence.
\end{proof}

\begin{lem} \label{EstWithP}
Let $i$ be a positive integer and suppose that $x\in X$ has the property that
$\|x\|\le C$ and $|x(\xi)|\le \delta$ whenever $\weight \xi=m_i$. Then for any
$s$ and any $\gamma$ of weight $m_i$ we have
$$
|\langle e^*_\gamma,P_{(s,\infty)}x\rangle|\le 2\delta + 3Cm_i^{-1}.
$$
\end{lem}
\begin{proof}
As before we consider the evaluation analysis
$$
e_\gamma^*= \sum_{r=1}^ad^*_{\xi_r} + m_i^{-1} \sum_{r=1}^a b^*_r\circ
P_{(p_{r-1},\infty)}.
$$

If $s\ge p_a$ then $P^*_{(s,\infty)}e^*_\gamma=0$. If $0<s<p_1$, by applying
$P^*_{(s,\infty)}$ to each of the terms in the evaluation analysis,  we see
that
$$
P^*_{(s,\infty)}e^*_\gamma = e^*_\gamma - m_i^{-1}P^*_{(0,s]}b^*_1,
$$
which leads to
$$
|\langle e_\gamma^*,P_{(s,\infty)}x_k\rangle| \le \delta
+m_i^{-1}\|b^*_1\|\,\|P_{(p_1,s]}\|\,\|x_k\|\le \delta + 3Cm_i^{-1},
$$
by our assumptions.

In the remaining case, there is some $t$ with  $1\le t<a$ such that $p_t\le s$
while $p_{t+1}>s$. We may rewrite the evaluation analysis of $\gamma$ as
$$
e_\gamma^*= e^*_{\xi_t} + \sum_{r=t+1}^ad^*_{\xi_r} + m_i^{-1} \sum_{r=t+1}^a
b^*_r\circ P_{(p_{r-1},\infty)},
$$
which gives us
$$
P^*_{(s,\infty)}e^*_\gamma= e^*_\gamma -
e^*_{\xi_t}-m_i^{-1}P^*_{(p_{t},s]}b^*_{t+1} .
$$
When we recall that $\weight \xi_t=\weight \gamma$ this yields
$$
|\langle e_\gamma^*, P_{(s,\infty)}x_k\rangle| \le 2\delta
+3Cm_i^{-1},
$$
as above.
\end{proof}

\begin{prop}[Basic Inequality]\label{BasicIneq}
Let $(x_k)_{k\in I}$ be a $C$-RIS, let $\lambda_k$ be real numbers, let $s$ be
a natural number and let $\gamma$ be an element of $\Gamma$.  There exist
$k_0\in I$ and and a functional $g^*\in W[(\mathscr A_{3n_j},m_j^{-1})_{j\in
\nats}]$ such that:
\begin{enumerate}
\item either $g^*=0$ or $\weight(g^*)=\weight(\gamma)$ and $\supp g^*\subseteq \{k\in I:k>k_0\}$ ;
 \item $\displaystyle |\langle e^*_\gamma,
P_{(s,\infty)}\sum_{k\in I}\lambda_kx_k\rangle|\le 5C|\lambda_{k_0}|
+ 5C\langle g^*,\sum_{k}|\lambda _k| e_k\rangle. $
\end{enumerate}
Moreover, if $j_0$ is such that
$$
|\langle e^*_{\xi},\sum_{k\in J}\lambda_k x_k\rangle|\le 2C\max_{k\in
J}|\lambda_k|,
$$
for all subintervals $J$ of $I$ and all $\xi\in \Gamma$ of weight $m_{j_0}$,
then we may choose $g^*$ to be in $W[(\mathscr A_{3n_j},m_j^{-1})_{j\ne j_0}]$.
\end{prop}
\begin{proof}
We proceed by induction of the rank of $\gamma$, noting that if $\gamma$ is  of
rank 1 we have $P^*_{(s,\infty)}e^*_\gamma=0$ whenever $s\ge 1$, so that
$$
\big\langle e^*_\gamma, P_{(s,\infty)}\sum_{k\in I}\lambda_kx_k\big\rangle =
\begin{cases} 0\qquad \qquad\text{if }r\ge 1\\
\lambda_1x_1(\gamma)\quad\text{if } r=0.
\end{cases}
$$
Thus $k_0=1$ and $g^*=0$ have the desired property.

Now consider an element $\gamma$ of rank greater than 1, of age $a$ and of
weight $m_h$.  Taking $(j_k)$ to be a sequence as in the definition of a RIS,
we shall suppose that there is some $l\in I$ such that $j_l\le h <j_{l+1}$.
(The cases where $h<j_k$ for all $k\in I$ and where $h\ge j_{k+1}$ for all
$k\in I$ are simpler.)

We split the summation over $k$ into three parts as follows:
$$
\langle e^*_\gamma, P_{(s,\infty)}\sum_{k\in
I}\lambda_kx_k\rangle=\sum_{I\ni k<l}\lambda_k\langle e^*_\gamma,
P_{(s,\infty)}x_k\rangle+\langle e^*_\gamma,
P_{(s,\infty)}\lambda_lx_l\rangle+\langle e^*_\gamma,
P_{(s,\infty)}\!\sum_{I\ni k>l}\!\!\lambda_kx_k\rangle
$$
and estimate the three terms separately.

When $k<l$ we have $h\ge j_l\ge j_{k+1}$ so that
$$
|\langle e^*_\gamma, P_{(s,\infty)}\lambda_kx_k\rangle|\le
3Cm_h^{-1}|\lambda_k|\le 3Cm_{j_k}^{-1}|\lambda_k|,
$$
by Lemma~\ref{RISEst}.  Thus
\begin{align*}
|\sum_{I\ni k<l}\lambda_k\langle e^*_\gamma, P_{(s,\infty)}x_k\rangle|\le
3C\sum_{k<l}m_{j_k}^{-1}|\lambda_k| \le 3C\sum_{j=1}^\infty
m_j^{-1}\max_{k<l}|\lambda_k|
 \le C\max_{k<l}|\lambda_k|.
\end{align*}
For the second term, we have the immediate estimate
$$
|\langle e^*_\gamma, P_{(s,\infty)}\lambda_lx_l\rangle|\le \|
P_{(s,\infty)}\|\,|\lambda_l|\,\|x_l\|\le 3C|\lambda_l|.
$$
Thus putting the first two terms together we have
\begin{equation}
|\langle e^*_\gamma, P_{(s,\infty)}\sum_{k\le l}\lambda_kx_k\rangle|\le
C\max_{k<l}|\lambda_k|+ 3C|\lambda_l|\le 4C|\lambda_{k_0}|,\label{BasIneq1}
\end{equation}
for a suitably chosen $k_0\le l$.

We now have to estimate the last term
 $$
  |\langle e^*_\gamma,\sum_{k\in I'}\lambda_kx'_k\rangle|,
 $$
 where $I'=\{k\in I:k>l\}$ and $x'_k = P_{(s,\infty)}x_k$. We shall use the evaluation analysis of $\gamma$
$$
e^*_\gamma = \sum_{r=1}^a d^*_{\xi_r} + m_h^{-1} \sum_{r=1}^a b^*_r\circ
P_{(p_{r-1},\infty)}.
$$

 Let $I'_0=\{k\in I':\ran x'_k\text{ contains }
\rank\xi_r\text{ for some } r\}$ noting first that $\#I'_0\le a$ and secondly
that for $k\in I'\setminus I'_0$ the interval $\ran x'_k$ meets
$(p_{r-1},p_r)$ for at most one value of $r$.  If we set $I'_r=\{k\in I': \ran
x_k\text{ meets } (p_{r-1},p_r]\text{  but no other }(p_{r'-1},p_{r'})\}$ then
each $I'_r$ is a subinterval of $I'$ and we have
\begin{align*}
\langle e^*_\gamma,x'_k\rangle &=  m_h^{-1}\langle
b^*_r,P_{(p_{r-1},\infty)}x'_k\rangle
=m_h^{-1}\langle b^*_r, P_{(s\vee p_{r-1},\infty)}x_k\rangle \\
\intertext{if $ k\in I'_r$, while} \langle e^*_\gamma,x'_k\rangle &=
    0 \qquad\qquad\text{if }k\in I'\setminus \bigcup_r I'_r
    \end{align*}
Thus $\displaystyle
 \langle e^*_\gamma,\sum_{k\in I'}\lambda_kx'_k\rangle= \langle e^*_\gamma,\sum_{k\in
 I'_0}\lambda_kx'_k\rangle
  +m_h^{-1}\sum_{r=1}^a \langle b^*_r, \sum_{k\in I'_r}\lambda _kx'_k\rangle.
$

Applying Lemma~\ref{RISEst}, we see that
\begin{equation}
|\langle e^*_\gamma,\sum_{k\in I'}\lambda_kx'_k\rangle|\le 5Cm_h^{-1}\sum_{k\in
I_0}|\lambda_k| + m_h^{-1}|\sum_{r=1}^a \langle b^*_r, \sum_{k\in I'_r}\lambda
_kx'_k\rangle|.\label{BasIneq2}
 \end{equation}
 Now, for each $r$, the
functional $b^*_r$ is a convex combination of functionals $\pm
e^*_\eta$ with $p_{r-1}<\rank \eta<p_r$, so we may choose $\eta_r$
to be such an $\eta$ with
$$
|\langle b^*_r,\sum_{k\in I'_r}\lambda _kx'_k\rangle|\le |\langle e^*_{\eta_r},
\sum_{k\in I'_r}\lambda _kx'_k\rangle|.
$$
For each $r$, we may apply our inductive hypothesis to the element $\eta_r\in
\Gamma$ and the RIS $(x_k)_{k\in I'_r}$, obtaining $k_r\in I'_r$ and $g^*_r\in
W[(\mathscr A_{3n_j},m_j^{-1})_{j\in \nats}]$ supported on $\{k\in I'_r:
k>k_r\}$ satisfying
\begin{equation}\label{BasIneq3}
|\langle e^*_{\eta_r}, P_{(s\vee p_r,\infty)}\sum_{k\in
I'_r}\lambda_kx_k\rangle|\le 5C|\lambda_{k_r}| + 5C\langle
g^*_r,\sum_{k\in I'_r}|\lambda _k| e_k\rangle.
\end{equation}
We now define $g^*$ by setting
 $$
g^* = m_h^{-1}(\sum_{k\in I'_0} e^*_{k} + \sum_{r=1}^a (e^*_{k_r}+g^*_r)).
$$
This is a sum, weighted by $m_h$, of at most $3n_h$ functionals in $W[(\mathscr
A_{3n_j},m_j^{-1})_{j\in \nats}]$, supported by disjoint intervals,  and is
hence itself in $W[(\mathscr A_{3n_j},m_j^{-1})_{j\in \nats}]$. Putting
together \ref{BasIneq1}, \ref{BasIneq2} and \ref{BasIneq3}, we finally obtain
\begin{align*}
|\langle e^*_\gamma, P_{(s,\infty)}(\sum_{k\in I}\lambda_kx_k\rangle|&\le
4C|\lambda_{k_0}|+ 5Cm_h^{-1}\sum_{k\in I'_0}|\lambda_k|
+m_h^{-1}|\sum_{r=1}^a \langle b^*_r,\sum_{k\in I'_r}\lambda _kx'_k\rangle|\\
&\le 4C|\lambda_{k_0}|+ 5Cm_h^{-1}\sum_{k\in I'_0}|\lambda_k|
+m_h^{-1}|\sum_{r=1}^a \langle e^*_{\eta_r}, P_{(s,\infty)} \sum_{k\in I'_r}\lambda _kx_k\rangle|\\
&\le 4C|\lambda_{k_0}|+ 5Cm_h^{-1}\biggl(\sum_{k\in I'_0}|\lambda_k|
+\sum_{r=1}^a ( |\lambda_{k_r}| + \langle g^*_r,\sum_{k\in
I'_r}|\lambda _k| e_k\rangle)\biggr)\\
&\le 5C|\lambda_{k_0}|+ 5C\langle g^*,\sum_{k\in I'}|\lambda_k|e_k\rangle.
\end{align*}

If $j_0$ satisfies the additional condition set out in the statement of the
theorem, we proceed by the same induction.  The base case certainly presents no
problem and if $\weight \gamma = m_h$ with $h=j_0$ we have a simple way to
estimate
$$
\langle e^*_\gamma, P_{(s,\infty)}\sum_{k\in I}\lambda _k x_k\rangle.
$$
Indeed there is at most one value of $k$, $l$ say, for which $s$ is in $\ran
x_k$ and $P_{(s,\infty)}x_k=0$ for $k<l$. If we set $J=\{k\in I: k>l\}$ we then
have
$$
|\langle e^*_\gamma,P_{(s,\infty)}\sum_{k\in I}\lambda _k
x_k\rangle|\le |\lambda_l|\|P_{(s,\infty)}\|\|x_l\|+
|e^*_\gamma(\sum_{k\in J}\lambda _k x_k)|,
$$
By our usual estimate $\|P_{(s,\infty)}\|\le 3$ and the assumed additional
condition, this is at most $5C|\lambda_{k_0}| $ for some $l\le k_0\in I$.  We
can then take $g^*=0$.
\end{proof}

\begin{cor}\label{RISDom} Any RIS is dominated by the unit vector
basis of $T[(\mathscr A_{3n_j},m_j^{-1})_{j\in \nats}]$. More precisely, if
$(x_k)$ is a $C$-RIS then, for any real $\lambda_k$, we have
$$
\|\sum_{k}\lambda_kx_k\|\le 10C\|\sum_k\lambda_ke_k\|,
$$
where the norm on the right hand side is taken in $T[(\mathscr
A_{3n_j},m_j^{-1})_{j\in \nats}]$.
\end{cor}

As well as this domination result, we shall need the following more precise
lemma.

\begin{prop}\label{RISUpperEst}
Let $(x_k)_{k=1}^{n_{j_0}}$ be a $C$-RIS.  Then
\begin{enumerate}
\item
For every $\gamma\in \Gamma$ with $\weight \gamma=m_h$ we have
$$
|n_{j_0}^{-1}\sum_{k=1}^{j_0} x_k(\gamma)|\le \begin{cases}
11Cm_{j_0}^{-1}m_h^{-1}\qquad\qquad\text{if }h<j_0\\
5Cn_{j_0}^{-1} + 5Cm_h^{-1} \qquad\text{ if } h\ge j_0 \end{cases}
$$
In particular,
$$
|n_{j_0}^{-1}\sum_{k=1}^{j_0} x_k(\gamma)|\le 6Cm_{j_0}^{-2},
$$
if $h>j_0$ and
$$
 \|n_{j_0}^{-1}\sum_{k=1}^{n_{j_0}}x_k\|\le 6Cm_{j_0}^{-1}.
$$
\item
If $\lambda_k$ ($1\le k\le n_{j_0}$) are scalars with $|\lambda_k|\le 1$ and
having the  property that
$$
|\sum_{k\in J} \lambda_kx_k(\gamma)|\le 2C\max_{k\in J}|\lambda_k|,
$$
for every $\gamma$ of weight $m_{j_0}$ and every interval
$J\subseteq\{1,2,\dots,n_{j_0}\}$, then
$$
 \|n_{j_0}^{-1}\sum_{k=1}^{j_0}\lambda_kx_k\|\le 6Cm_{j_0}^{-2}.
$$
\end{enumerate}
\end{prop}
\begin{proof}
This is a direct application of the Basic Inequality, with all the coefficients
$\lambda_k$ equal to $n_{j_0}^{-1}$. Indeed, for (1) there exists $g^*\in
W[(\mathscr A_{3n_j},m_j^{-1})_{k\in \nats}]$ (either zero or of weight $m_h$)
such that
$$
|n_{j_0}^{-1}\sum_{k=1}^{n_{j_0}} x_k(\gamma)|\le 5Cn_{j_0}^{-1} +
5Cg^*(n_{j_0}^{-1}\sum_{k=1}^{n_{j_0}} e_k)).
$$
Using Lemma~\ref{MTLem} to estimate the term involving $g^*$, we obtain
$$
|n_{j_0}^{-1}\sum_{k=1}^{n_{j_0}} x_k(\gamma)|\le
\begin{cases} 5Cn_{j_0}^{-1}+10Cm_{j_0}^{-1}m_h^{-1}\quad \text{if }h<j_0\\
              5Cn_{j_0}^{-1} + 5Cm_h^{-1}\qquad\qquad\text{if } h\ge j_0.
\end{cases}
 $$
The formulae given in (1) follow easily when we note that $n_{j_0}$ is (much)
larger than $5m_{j_0}^2$ when $j_0\ge 2$.

If the scalars $\lambda_k$ satisfy the additional condition, then the $g^*$
whose existence is guaranteed by the Basic Inequality may be taken to be in
$W[(\mathscr A_{3n_j},m_j^{-1})_{j\ne j_0}]$ so that the second part of
Lemma~\ref{MTLem} may be applied, yielding
$$
|n_{j_0}^{-1}\sum_{k=1}^{n_{j_0}} x_k(\gamma)|\le
\begin{cases} 5Cn_{j_0}^{-1}+10Cm_{j_0}^{-2}m_h^{-1}\quad \text{if }h<j_0\\
              5Cn_{j_0}^{-1} + 5Cm_h^{-1}\qquad\qquad\text{if } h> j_0.
\end{cases}
 $$
This leads easily to the claimed estimate for
$\|n_{n_{j_0}}^{-1}\sum_{k=1}^{j_0}\lambda_kx_k\|$.
\end{proof}

It turns out that in our space there are three useful types of RIS. One of
these is based on an idea that will be familiar from other constructions, that
of introducing long $\ell_1$-averages. We defer our discussion of this
construction until the next section. We shall deal first with the other two
types of RIS, which involve the $\mathscr L_\infty$ structure of our space, and
provide the extra tool that we eventually use to solve the scalar-plus-compact
problem.

We have already remarked that the {\em support} of an element of $X$ is not of
great interest --- indeed the support of any nonzero element of $X$ is an
infinite set, and contains elements $\gamma$ of $\Gamma$ of all possible
weights.  There is, however, a related notion which is of much use. Recall that
an element $x$ whose range is contained in the interval $(p,q]$ can be
expressed as $i_q(u)$ where $u\in \ell_\infty(\Gamma_q)$ and
$\supp(u)\subseteq\Gamma_q\setminus \Gamma_p$. It turns out that the support of
$u$ contains a lot of information about $x$.  We shall refer to $\supp(u)$ as
the {\em local support}. A formal (and unambiguous) definition may be
formulated as follows.

\begin{defn}
Let $x$ be an element of $\bigoplus_nM_n$ and let $q=\max\ran x$; thus $x$ may
be expressed as $i_q(u$)with $u=x\restriction \Gamma_q$. The subset $\supp
u=\{\gamma\in \Gamma_q:x(\gamma)\ne 0\}$ is defined to be the {\em local
support} of $x$.
\end{defn}

The following easy lemma uses an idea that has already occurred in
Lemma~\ref{RISEst}.

\begin{lem}\label{LocSuppWeight}
Let $\gamma\in \Gamma$ be of weight $m_h$ and assume that $\weight(\xi) \ne
m_h$ for all $\xi$ in the local support of $x$. Then $|x(\gamma)|\le
3m_h^{-1}\|x\|.$
\end{lem}
\begin{proof}
Let $q=\max\ran x$ so that $x=i_q(x\restriction \Gamma_q)$ and, by hypothesis,
$\weight \xi\ne m_h$ whenever $\xi\in \Gamma_q$ and $x(\xi)\ne 0$. If
$\rank\gamma\le q$ we thus have $x(\gamma)=0$ and there is nothing to prove.
Otherwise we consider the evaluation analysis of $\gamma$
$$
e^*_\gamma = \sum_{r=1}^a d^*_{\xi_r} + m_h^{-1}\sum_{r=1}^a b^*_r\circ
P_{(p_{r-1},\infty)}
$$
and let $s$ be chosen maximal subject to $p_s=\rank \xi_s\le q$. (Since
$\gamma=\xi_a$ such an $s$ certainly exists.)  For $r\ge s$ we have $r>\max\ran
x$, whence $d^*_{\xi_r}(x) = 0$ and $P_{(p_r,\infty)}x=0$. Thus
$$
x(\gamma)= \langle e^*_{\gamma},x\rangle = \begin{cases}
m_h^{-1}\langle b^*_s,P_{(p_{s-1},\infty)}x\rangle+\langle
e^*_{\xi_{s-1}},x\rangle\quad\text{if
}s>1\\
m_h^{-1}\langle b^*_1,x\rangle = m_h^{-1}\langle
b^*_1,P_{(p_0,\infty)}x\rangle\quad\text{ if }s=1.
\end{cases}
$$
Since, in the first of the above cases, we have $\rank \xi_{s-1}<q$ and
$\weight \xi_{s-1}=m_h$, which imply $e^*_{\xi_{s-1}}(x)=0$, we deduce that in
both cases
$$
|x(\gamma)| =m_h^{-1}|\langle b^*_s,P_{(p_{s-1},\infty)}x\rangle|\le
3m_h^{-1}\|x\|.
$$
\end{proof}

We can now introduce two classes of block sequence, characterized by the
weights of the elements of the local support.

\begin{defn}
We shall say that a block sequence $(x_k)_{k\in \nats}$ in $X$ has {\em bounded
local weight} if there exists some $j_1$ such that $\weight \gamma\le m_{j_1}$
for all $\gamma$ in the local support of $x_k$, and all values of $k$.  We
shall say that $(x_k)_{k\in\nats}$ has {\em rapidly increasing local weight}
if, for each $k$ and each $\gamma$ in the local support of $x_{k+1}$, we have
$\weight \gamma>m_{i_k}$ where $i_k=\max\ran x_k$.
\end{defn}

\begin{prop}\label{LocWeightRIS}
Let $(x_k)_{k\in \nats}$ be a bounded block sequence.  If either $(x_k)$ has
bounded local weight, or $(x_k)$ has rapidly increasing local weight, the
sequence $(x_k)$ is a RIS.
\end{prop}
\begin{proof}
We start with the case of rapidly increasing local weight and let $m_{j_k}$ be
the minimum weight of an element $\gamma$ in the local support of $x_k$.  By
hypothesis, $j_{k+1} > \max\supp x_k$ so that RIS condition (2) is satisfied.
Also, if $h<j_k$ and $\gamma$ is of weight $m_h$ then $|x_k(\gamma)|\le
3m_h^{-1}\|x_k\|$ by Lemma~\ref{LocSuppWeight}.  So $(x_k)$ is a $C$-RIS with
$C=3\sup_k\|x_k\|.$

Now let us suppose that $\weight \gamma\le m_{j_1}$ for all $\gamma$ in the
local support of $x_k$ and all $k$.  For $k\ge 2$ define $j_k=1+\max \supp
x_{k-1}$, thus ensuring that RIS condition (2) is satisfied.  If $\weight
\gamma=m_h$ where $h<j_k$ there are two possibilities: if $i>j_1$ then
$|x_k(\gamma)|\le 3m_i^{-1}\|x_k\|$ by Lemma~\ref{LocSuppWeight}; if $i\le j_1$
then $|x_k(\gamma)|\le \|x_k\|\le m_i^{-1}m_{j_1}\|x_k\|$.  Thus $(x_k)$ is a
$C$-RIS, where $C$ is the (possibly quite large) constant
$m_{j_1}^{-1}\sup_k\|x_k\|$.
\end{proof}

\begin{prop}\label{RIStoBlock}
Let $Y$ be any Banach space and $T:X(\Gamma)\to Y$ be a bounded linear
operator.  If $\|T(x_k)\|\to0$ for every RIS $(x_k)_{k\in \nats}$ in
$X(\Gamma)$  then $\|T(x_k)\|\to 0$ for every bounded block sequence sequence
in $X(\Gamma)$.
\end{prop}
\begin{proof}
It is enough to consider a bounded block sequence $(x_k)$ and show that there
is a subsequence $(x'_j)$ such that $\|T(x'_j)\|\to 0$. We may write
$x_k=i_{q_k}(u_k)$ with $u_k=x_k\restriction \Gamma_{q_k}$ supported by
$\Gamma_{q_k}\setminus \Gamma_{q_{k-1}}$. For each $k$ and each $N\in \nats$,
we split $u_k$ as $v^N_k+w^N_k$, where, for $\gamma\in \Gamma_{q_k}$,
 \begin{align*}
v^N_k(\gamma) &=\begin{cases} u_k(\gamma)\text{ if }\weight
\gamma\le m_N\\
0\quad\qquad\text{otherwise}\end{cases}\\
w^N_k(\gamma) &=\begin{cases} u_k(\gamma)\text{ if }\weight \gamma> m_N\\
0\quad\qquad\text{otherwise}\end{cases}
\end{align*}
and set
$$
\quad y^N_k = i_{q_k}(v^N_k),\qquad z^N_k = i_{q_k}(v^N_k).
$$
We notice that $\|y_k^N\|\le \frac32\|v_k^N\|\le \frac32\|x_k\|$, with a
similar estimate for $\|z^N_k\|$, so that the sequences $(y^N_k)_k$ and
$(z^N_k)_k$ are bounded. We note also that $\weight\gamma\le N$ for all $\gamma
$ in the local support of $y^N_k$ and $\weight\gamma> N$ for all $\gamma $ in
the local support of $z^N_k$

So for each $N$, the sequence $(y^N_k)$ has bounded local weight and is thus a
RIS, by Proposition~\ref{LocWeightRIS}. By hypothesis, $\|T(y^N_k)\|\to 0$ for
each $N$.   Hence we can choose a sequence $(k_n)$ tending to $\infty$ such
that $\|T(y^n_{k_n}\|\to 0$.  If we put $n_1=1$ and then, recursively, set
$n_{j+1}= q_{k_{n_j}}$, it is easy to see that the sequence
$(z^{n_j}_{k_{n_j}})$ has rapidly increasing local weight.  Thus this sequence
is a RIS and we hence have $\|T(z^{n_j}_{k_{n_j}})\|\to 0$.  Since
$x_{k_{n_j}}=y^{n_j}_{k_{n_j}}+z^{n_j}_{k_{n_j}}$ , we have found a subsequence
$(x'_j)=(x_{k_{n_j}})$ of $(x_k)$ with $\|T(x'_j)\|\to 0$.
\end{proof}

The above proposition will play an important role in proving compactness of
operators, but in the mean time we shall use it to give our promised proof that
the dual of $X$ is $\ell_1$.  There is an alternative approach using
$\ell_1$-averages.

\begin{prop}
The dual of $X(\Gamma)$ is $\ell_1(\Gamma)$.
\end{prop}
\begin{proof}
As we have already noted in Theorem~\ref{BDThmBis} it is enough to show that
the FDD $(M_n)$ is shrinking, that is to say, that every bounded block sequence
in $X$ is weakly null.  So let $\phi$ be an element of $X^*$.  By the upper
estimate of Proposition~\ref{RISUpperEst} we see that $\phi(x_k)\to 0$ for
every RIS $(x_k)_{k\in \nats}$. Now Proposition~\ref{RIStoBlock}, applied with
$T=\phi$, shows that $\phi(x_k)\to 0$ for every bounded block sequence $(x_k)$.
\end{proof}

\begin{remark}\label{ReflRem} We can see that $X(\Gamma)$ has many reflexive subspaces.
Indeed, suppose that $(q_n)$ is an increasing sequence of natural numbers, and
that, for each $n$, $F_n$ is a finite dimensional subspace of
$\bigoplus_{q_n<k<q_{n+1}}M_k$. Then $(F_n)$ is an FDD for the subspace
$W=\overline{\bigoplus_{n\in\nats} F_n}$, and $(F_n)$ is shrinking because
$(M_n)$ is.  But $(F_n)$ is also boundedly complete, by the lower estimate of
Proposition~\ref{LowerEst}. Thus $W$ is reflexive.  We shall see later in
section~\ref{HIDuals} that $W^*$ is hereditarily indecomposable whenever $W$ is
a subspace of this type in the space $\XK$.
\end{remark}

\section{$\ell_1$-averages, Exact Pairs and the HI
property}\label{ellOneExact}

In the first part of this section, we shall still only be using the assumptions
\ref{DeltaUpperAssump} and \ref{DeltaLowerAssump}, so that our results will
apply when $X$ is either of the spaces $\BmT$ and $\XK$.  The special
properties of the second of these spaces will come into play only from
Definition~\ref{DepSeq} onwards.

\begin{defn}\label{ellOneDef}
An element $x$ of $X$ will be called a {\em $C$-$\ell_1^n$ average}
if there exists a block sequence $(x_i)_{k=1}^n$ in $X$ such that
$x=n^{-1}\sum_{k=1}^n x_k$ and $\|x_k\|\le C$ for all $k$. We say
that $x$ is a {\em normalized }$C$-$\ell_1^n$ average if, in
addition, $\|x\|=1$.
\end{defn}

A standard argument (c.f. II.22 of \cite{AT}) using the lower estimate of
Lemma~\ref{LowerEst} and Lemma~\ref{ellOneAve} leads to the following.

\begin{lem}\label{ExistAve}
Let $Z$ be any block subspace of $X$.  For any $n$ and and $C>1$, $Z$ contains
a normalized $C$-$\ell_1^n$ average.
\end{lem}
\begin{proof}
Write $C=(1-\epsilon)^{-1}$ and choose an integer $l$ with
$n(1-\epsilon/n)^l<1$; next choose $j$ sufficiently large as to ensure that
$n_{2j}>(2m_{2j})^l$; finally let $k$ be minimal subject to
 $$
m_{2j}<(1-\epsilon/n)^{-k}
 $$
 Since $\frac12(1-\epsilon/n)^{-k}\le(1-\epsilon/n)^{-k+1}\le m_{2j}$ we have
$$
n_{2j}>(2m_{2j})l\ge (1-\epsilon/n)^{-kl}>n^k.
$$

If $(x_i)$ is any normalized skipped-block sequence in $Z$, we can apply
Lemma~\ref{LowerEst} to see that
$$
\|\sum_{i=1}^{n^k} x_i\|\ge m_{2j}^{-1}n^k>(n-\epsilon)^k.
$$
It now follows from Lemma~\ref{ellOneAve} that there are normalized successive
linear combinations $y_1,\dots,y_n$ of $(x_i)$ such that
$$
\|\sum_{i=1}^n a_iy_i\|\ge (1-\epsilon)\sum_{i=1}^n|a_i|,
$$
for all real $a_i$.  In particular, there is a normalized $C$-$\ell_1^n$
average.
\end{proof}

\begin{lem}\label{SmallWeightLongAve}
Let $x$ be a $C$-$\ell_1^{n_j}$ average.  For all $\gamma\in \Gamma$ we have
$|\langle d^*_\gamma,x\rangle|\le 3Cn_j^{-1}$.  If $\gamma$ is of weight $m_i$
with $i<j$ and $p\in \nats$ then $|x(\gamma)|\le 2Cm_i^{-1}$.
\end{lem}
\begin{proof}
Let $x= n_j^{-1}\sum_{k=1}^{n_j}x_k$, as in the definition of a $C$-$\ell_1^n$
average. For any $\gamma$ there is some $k$ such that $\langle
d^*_{\gamma},x\rangle=n_j^{-1}\langle d^*_{\gamma},x_k\rangle$. Thus
$$
|\langle d^*_\gamma,x\rangle| \le n_j^{-1}\|d^*_{\gamma}\|\, \|x_k\|\le
3Cn_j^{-1}.
$$

Let us now consider the case where $\weight \gamma=m_i$, with $i<j$. From the
evaluation analysis
$$
e^*_\gamma = \sum_{r=1}^a d^*_{\xi_r} + m_i^{-1} \sum_{r=1 }^a b^*_{r}\circ
P_{(p_{r-1},\infty)},
$$
it follows that
\begin{equation}\label{SWLA1}
|x(\gamma)| \le \sum_{r=1}^a |\langle d^*_{\xi_r},x\rangle| + m_i^{-1}
\sum_{r=1}^a \|P_{(p_{r-1},p_r)}x\|.
\end{equation}

By what we have already observed, we have
\begin{equation}\label{SWLA2}
 \sum_{r=1}^a |\langle d^*_{\xi_r},x\rangle| \le 3Can_j^{-1}.
 \end{equation}

To estimate the second term in \eqref{SWLA1} we follow the argument of page 33
of \cite{AT}, letting $I_r$ (resp. $J_r$) be the set of $k$ such that $\ran
x_k$ is contained in (resp. meets) the interval $(p_{r-1},p_r)$.  We have
$\#J_r\le \#I_r+2$ and $\sum_r \#I_r\le n_j$. Moreover, for each $r$, we have
$P_{(p_{r-1},p_r)}x_k=x_k$ if $k\in I_r$, while $P_{(p_{r-1},p_r)}x_k=0$ if
$k\notin J_r$ and
$$
\|P_{(p_{r-1},p_r)}x_k\|\le 4\|x_k\|\le 4C\quad\text{if } k\in J_r\setminus
I_r.
$$
It follows that $$
 \|P_{(p_{r-1},p_{r}]}x\|\le n_j^{-1} (C\#I_r +
8C))\le Cn_j^{-1} (\#I_r +8).
$$
Summing over $r$ leads us to
\begin{equation}\label{SWLA3}
 \sum_{r\le a}
\|P_{(p_{r-1},p_{r}]}x\|\le Cn_j^{-1}(n_j + 8a).
\end{equation}
Combining our inequalities, and using the fact that $a\le n_i$ we obtain
$$
|x(\gamma)| \le 3Can_j^{-1}+ m_i^{-1}n_j^{-1}(Cn_j + 8Ca) \le
Cm_i^{-1}+5Cn_in_j^{-1}<2Cm_i^{-1}.
$$
\end{proof}

\begin{lem}\label{AveRIS}
Let $I$ be an interval in $\mathbb N$, et $(x_k)_{k\in I}$ be a block sequence
in $X$ and let $(j_k)_{k\ge 1}$ be an increasing sequence of natural numbers.
Suppose that, for each $k$, $x_k$ is a $C$-$\ell_1^{n_{j_k}}$-average and that
$j_{k+1} > \max \ran x_k$. Then $(x_k)$ is a $2C$-RIS.
\end{lem}
\begin{proof}
We just have to prove RIS condition (3) and this is an immediate consequence of
Lemma~\ref{SmallWeightLongAve}.
\end{proof}

\begin{cor}\label{ExistRIS}
Let $Z$ be a block subspace of $X$, and let $C>2$ be a real number. Then $X$
contains a normalized $C$-RIS.
\end{cor}
\begin{proof}
This is immediate from Lemmas~\ref{ExistAve} and \ref{AveRIS}.
\end{proof}

\begin{defn}
Let $C>0$ and let $\varepsilon\in \{0,1\}$.  A pair $(x,\gamma)\in X\times
\Gamma$ is said to be a $(C,j,\varepsilon)$-{\em exact pair} if
 \begin{enumerate}
 \item $|\langle d^*_\xi,x\rangle|\le Cm_{j}^{-1}$ for all $\xi\in \Gamma$;
 \item $\weight \gamma = m_{j}$,
$\|x\|\le C$, $x(\gamma) = \varepsilon$;
 \item for every element $\gamma'$ of
$\Gamma$ with $\weight\gamma'= m_{i}\ne m_j$, we have
$$
|x(\gamma')| \le \begin{cases} Cm_i^{-1} \qquad\text{if } i<j\\
                 Cm_{j}^{-1} \qquad\text{if }i>j.
                            \end{cases}
$$
\end{enumerate}
\end{defn}

\begin{remark}
It is an immediate consequence of Lemma~\ref{EstWithP} that a
$(C,j,\varepsilon)$ exact pair also satisfies the estimates
$$
|\langle e^*_\gamma,P_{(s,\infty)}x\rangle| \le \begin{cases} 5Cm_i^{-1} \qquad\text{if } i<j\\
                 5Cm_{j}^{-1} \qquad\text{if }i>j
                            \end{cases}
$$
for elements $\gamma $ of $\Gamma$ with $\weight \gamma=m_i\ne m_j$.
\end{remark}

It will be seen that these estimates, as well as those in the definition, have
much in common with those of Lemma~\ref{RISEst}. Our first task is to show how
we can construct $(C,2j,1)$-exact pairs, starting from a RIS.

\begin{lem}\label{RISOneExact}
Let $j$ be a positive integer and let $(x_k)_{k=1}^{n_{2j}}$ be a
skipped-block $C$-RIS, such that $\min\ran x_2\ge 2j$ and
$\|x_k\|\ge 1$ for all $k$.  Then there exists $\theta\in \reals$,
with $|\theta|\le 2$, and there exists $\gamma\in \Gamma$, such that
$(x,\gamma)$ is a $(22C,2j,1)$-exact pair, where $x$ is the weighted
sum
$$
 x=\theta m_{2j}n_{2j}^{-1}\sum_{k=1}^{n_{2j}}x_k.
$$
\end{lem}
\begin{proof}
We may apply the construction of Lemma~\ref{LowerEst} to obtain an element
$\gamma$ of $\Gamma$ of weight $m_{2j}$ such that
$$
n_{2j}^{-1}\sum_{k=1}^{n_{2j}} x_k(\gamma) \ge \half m_{2j}^{-1}.
$$
For a suitably chosen $\theta\in \reals$ with $0<\theta\le 2$ we have
$x(\gamma)=1$, where $x=\theta m_{2j} n_{2j}^{-1}\sum_{k=1}^{n_{2j}} x_k$. We
thus have condition (2) in the definition of an exact pair.

There is no problem establishing condition (1) since, for any $\xi$, there is
some $k$ satisfying $\langle d^*_\xi,x\rangle=\theta m_{2j}n_{2j}^{-1}\langle
d^*_\xi,x_k\rangle$. By RIS condition (1), $\|x_k\|\le C$ and we know that
$\|d^*_\xi\|\le 3$.  Hence $|\langle d^*_\xi,x_k\rangle| \le
6Cm_{2j}n_{2j}^{-1}<Cm_{2j}^{-1}.$

To establish condition (3), we shall use the fact that $(x_k)$ is a $C$-RIS and
apply Proposition~\ref{RISUpperEst}, with $j_0=2j$.  If $\weight \gamma'=m_i$
with $i\ne2j$, we thus have
 $$
 |x(\gamma)|=|\theta |m_{2j}n_{2j}^{-1}\sum_{k=0}^{n_{2j}} x_k(\gamma')\le \begin{cases}
 22 Cm_i^{-1}\qquad\qquad\qquad\qquad\qquad\qquad\qquad\text{if } i<2j\\
 10Cm_{2j}n_{2j}^{-1}+10Cm_{2j}m_i^{-1}<11Cm_{2j}^{-1}\quad\text{if } i>2j.
 \end{cases}
 $$

 \end{proof}

Using Lemma~\ref{ExistRIS} we now immediately obtain the following.

\begin{lem}\label{ExistExact}
If $Z$ is a block subspace of $X$ then  for every $j\in \nats$ there exists a
$(45,2j,1)$-exact pair $(x,\eta)$ with $x\in Z$.
\end{lem}

The proof of the following lemma, which we shall need in
Section~\ref{Operators}, is very similar.

\begin{lem}\label{RISZeroExact}
Let $(x_k)_{k=1}^{n_{2j}}$ be a skipped-block $C$-RIS, and let
$q_0<q_1<Q_2<\cdots<q_{n_{2j}}$ be  natural numbers such that $\ran
x_k\subseteq (q_{k-1},q_k)$ for all $k$. Let $z$ denote the weighted sum
$x=m_{2j}n_{2j}^{-1}\sum_{k=1}^{n_{2j}}x_k$. For each $k$ let $b^*_k$ be an
element of $B_{q_k-1,q_{k-1}}$ with $b_k(x_k)=0$. Then there exist $\zeta_i \in
\Delta_{q_i}$ ($1\le i\le n_{{2j}}$) such that the element
$\eta=\zeta_{n_{2j}}$ has analysis $(q_i,b^*_i,\zeta_i)_{1\le i\le n_{2j}}$ and
$ (z,\eta)$ $(12C,n_{2j},0)$-exact pair.
\end{lem}

We are finally ready to make use of the special conditions governing
``odd-weight'' elements of $\Gamma$.  We need to consider a special type of
rapidly increasing sequence whose members belong to exact pairs.

\begin{defn}\label{DepSeq}Consider the space $\XK=X(\Gamma)$ where
$\Gamma=\Gamma_{\mathrm K}$ as defined in \ref{XKDef}. We shall say that a
sequence $(x_i)_{i\le n_{2j_0-1}}$ is a $(C, 2j_0-1,\varepsilon)$-{\em
dependent sequence} if there exist $0=p_0<p_1<p_2<\cdots<p_{n_{2j-1}}$,
together with $\eta_i\in\Gamma_{p_i-1}\setminus \Gamma_{p_{i-1}}$ and $\xi_i\in
\Delta_{p_i}$ ($1\le i\le n_{2j_0-1}$) such that
\begin{enumerate}
\item for each $k$, $\ran x_k\subseteq (p_{k-1},p_k)$;
 \item the
element $\xi=\xi_{2j_0-1}$ of $\Delta_{p_{2j_0-1}}$ has weight $m_{2j_0-1}$ and
analysis $(p_i,e^*_{\eta_i},\xi_i)_{i=1}^{2n_{j_0}-1}$;
\item $(x_1,\eta_1)$ is a $(C,4j_1-2,\varepsilon)$-exact pair; \item
for each $2\le i\le n_{2j-1}$, $(x_i,\eta_i)$ is a $(C,4j_i,\varepsilon)$-exact
pair, with $\ran x_i\subseteq (p_{i-1},p_i)$.
\end{enumerate}
We notice that, because of the special odd-weight conditions in
\ref{XKDef}, we necessarily have $m_{4j_1-2} = \weight \eta_1
>n_{2j_0-1}^2$, and $\weight \eta_{i+1} =m_{4j_{i+1}}$, where
$j_{i+1}=\sigma(\xi_i)$ for $1\le i<n_{2j_0-1}$.
\end{defn}

\begin{lem}
A $(C, 2j_0-1,\varepsilon)$-{dependent sequence} in $\XK$ is a $C$-RIS.
\end{lem}
\begin{proof}
For each $i\ge 1$ we have $\max \ran x_i <p_i$ and $j_{i+1}=\sigma(\xi_i)>\rank
\xi_i=p_i$.  This establishes Condition (2) in the definition of a RIS.
Condition (3) follows from the the definition of a $C$-exact pair.
\end{proof}

\begin{lem}\label{PlusMinus}
Let $(x_i)_{i\le n_{2j_0-1}}$ be a $(C, 2j_0-1,1)$-{dependent sequence} in
$\XK$ and let $J$ be a sub-interval of $[1,n_{2j_0-1}]$. For any $\gamma'\in
\Gamma$ of weight $m_{2j_0-1}$ we have
$$
|\sum_{i\in I} (-1)^i x_i(\gamma')| \le 4C.
$$
\end{lem}
\begin{proof}
Let $\xi_i,\eta_i,p_i,j_i$ be as in the definition of a dependent sequence and
let $\gamma$ denote $\xi_{2j_0-1}$, an element of weight $m_{4j_0-1}$. Let
$(p_i',e^*_{\eta'_i},\xi'_i)_{1\le i\le a'}$) be the analysis  of $\gamma'$ and
let the weight of $\xi'_i$ be $m_{4j'_1-2}$ when $i=1$, $m_{4j'_i}$ when
$1<i\le a'$. We note that $a'\le n_{2j_0-1}$ because $\gamma'$ is of weight
$m_{2j_0-1}$. We may thus apply the tree-like property of Lemma~\ref{Treelike}
deducing that there exists $1\le l\le a'$ such that
$(p_i',\eta_i',\xi'_i)=(p_i,\eta_i,\xi_i)$ for $i<l$ while $j_k\ne j'_i$ for
all $l<i\le a'$ and all $1\le k\le {n_{2j_0-1}}$. Since
$$
e^*_{\gamma} \circ P_{(0,p_{l-1}]} =
e^*_{\xi_{l-1}}=e^*_{\xi_{l-1}}=e^*_{\gamma} \circ P_{(0,p_l)},
$$
we have
$$
 x_k(\gamma')=x_k(\gamma)=m_{2j_0-1}^{-1}e^*_{\eta_k}\circ
P_{p_{k-1},\infty)}x_k=m_{2j_0-1}^{-1}x_k(\eta_k)=m_{2j_0-1}^{-1},
 $$ for $1\le k<l$.

We may now estimate as follows
\begin{align*}
|\sum_{k\in I} (-1)^k x_k(\gamma')| &\le |\sum_{k\in J,\ k<l}m_{2j_0-1}^{-1}
(-1)^k x_k(\gamma')| + |x_l(\gamma')| +
\sum_{k\in J,\ k>l}|x_i(\gamma')|\\
&\le m_{2j_0-1}^{-1}|\sum_{k\in I,\ k<l} (-1)^k|+\|x_l\|+\sum_{k\in J,\
k>l}\sum_{i\le a'}|d^*_{\xi_i}x_k +
m_{2j_0-1}^{-1}e^*_{\eta'_i}\circ P_{(p'_{i-1},\infty)}x_k|\\
&\le 1+C+n_{2j_0-1}^2\max_{l<k\in J,\ i\le a'}|d^*_{\xi_i}x_k +
m_{2j_0-1}^{-1}e^*_{\eta'_i}\circ P_{(p'_{i-1},\infty)}x_k|.
\end{align*}

 Now we know that, provided $k>l$, $\weight \eta'_k
\ne \weight \eta_i$ for all $i$, so by the definition of an exact pair, we have
 \begin{align*}
 |d^*_{\xi'_k}(x_i)+m_{2j_0-1}^{-1}P_{(p_{k-1},\infty]}x_i(\eta'_k)|&\le
 C(\weight \eta_i)^{-1} + 5Cm_{2j_0-1}^{-1}\max\{(\weight \eta'_{k})^{-1}, (\weight \eta_{i})^{-1}\}\\
 &\le 2C \max\{(\weight \eta_1)^{-1},\weight \eta'_1)^{-1}\}\\
 &= 2C\max\{m_{4j_1-2}^{-1},m_{4j'_1-2}^{-1}\}\le 2C
 n_{2j_0-1}^{-2},
 \end{align*}
using the fact that $m_{4j_1-2}$ and $m_{4j'_1-2}$ are both at least
$ n_{2j_0-1}^2$. We now deduce the inequality $|\sum_{i\in J} (-1)^i
x_i(\gamma')| \le 4C$ as required.
\end{proof}

\begin{lem}\label{OneDepSeqEst}
Let $(x_i)_{i\le n_{2j-1}}$ be a $(C, 2j_0-1,1)$-{dependent sequence} in $\XK$.
Then $$
 \|n_{2j_0-1}^{-1}\sum_{i=1}^{n_{2j_0-1}}
x_i\|\ge m_{2j_0-1}^{-1}\quad \text{ but }\quad
\|n_{2j_0-1}^{-1}\sum_{i=1}^{n_{2j_0-1}} (-1)^ix_i\| \le 12Cm_{2j_0-1}^{-2}.
 $$
\end{lem}
\begin{proof}
Using the notation of Definition~\ref{DepSeq} is easy to show by induction on
$a$, as in Lemma~\ref{LowerEst}, that
$$
\sum_{i=1}^{a} x_i(\xi_a) = m_{2j_0-1}^{-1}a,
$$
whence we immediately obtain
$$
\|n_{2j_0-1}^{-1}\sum_{i=1}^{n_{2j_0-1}} x_i\|\ge\sum_{i=1}^{a}
x_i(\xi_{2j_0-1}) \ge m_{2j_0-1}^{-1}.
$$

To estimate $\|n_{2j_0-1}^{-1}\sum_{i=1}^{n_{2j_0-1}} (-1)^ix_i\|$ we consider
any $\gamma\in \Gamma$ and apply the second part of Lemma~\ref{RISUpperEst},
with $\lambda_i = (-1)^nn_{2j_0-1}^{-1}$ and with $2j_0-1$ playing the role of
$j_0$. Lemma~\ref{PlusMinus} shows that the extra hypothesis of the second part
of Lemma~\ref{RISUpperEst} is indeed satisfied, provided we replace $C$ by
$2C$. We deduce that $\|n_{2j_0-1}^{-1}\sum_{i=1}^{n_{2j_0-1}} (-1)^ix_i\| \le
12Cm_{2j_0-1}^{-2},$ as claimed.
\end{proof}

A very similar proof yields the following estimate, which we shall use in the
next section.

\begin{lem}\label{ZeroDepSeqEst}
Let $(x_i)_{i\le n_{2j-1}}$ be a $(C, 2j_0-1,0)$-{dependent sequence} in $\XK$.
Then $$ \|n_{2j_0-1}^{-1}\sum_{i=1}^{n_{2j_0-1}} x_i\| \le 4Cm_{2j_0-1}^{-2}.
 $$
\end{lem}

In the mean time, we may finish the proof of one of our main theorems.

\begin{lem}\label{HILem}
Let $Y$ and $Z$ be block subspaces of $\XK$.  Then, for each $\epsilon>0$,
there exist $y\in Y$ and $z\in Z$ with $\|y-z\|<\epsilon \|y+z\|$.
\end{lem}
\noindent{\em Proof.} We start by choosing $j_0,j_1$ with
$m_{2j_0-1}>540\epsilon^{-1}$ and $m_{4j_1-2}>n_{2j_0-1}^2$. Next we use
Lemma~\ref{ExistExact} to choose a $(45,m_{4j_1-2},1)$-exact pair
$(x_1,\eta_1)$ with $x_1\in Y$. Now, for some $p_1> \rank\eta_1\vee \max\ran
x_1$, we define $\xi_1\in \Delta_{p_1}$ to be $(p_1,m_{2j_0-1},e^*_{\eta_1})$.

 We now set
$j_2=\sigma(\xi_1)$ and choose a $(45,m_{4j_2},1)$-exact pair $(x_2,\eta_2)$
with $x_2\in Z$ and $\min\ran x_2>p_1$.  We pick $p_2> \rank \eta_2\vee\max\ran
x_2$ and take $\xi_2$ to be the element $(p_2,\xi_1,m_{2j_0-1}e^*_{\eta_2})$ of
$\Delta_{p_2}$. Notice that this tuple is indeed in $\Delta_{q_2+1}$ because we
have ensured that $\weight\eta_2=m_{4\sigma(\xi_1)}$.

Continuing in this way, we obtain a $(45,2j_0-1)$-dependent sequence $(x_i)$
such that $x_i\in Y$ when $i$ is odd and $x_i\in Z$ when $i$ is even.  We
define $y= \sum_{i \text{ odd}} x_i$ and $z= \sum_{i \text{ even}} x_i$, and
observe that, by Lemma~\ref{DepSeq},
\begin{align*}
 \|y+z\|&=\|\sum_{i=1}^{n_{2j_0-1}}x_i\|\ge n_{2j_0-1}m_{2j_0-1}^{-1},\quad\text{while}\\
\|y-z\|&=\|\sum_{i=1}^{n_{2j_0-1}}(-1)^ix_i\|\le 12\times45
n_{2j_0-1}m_{2j_0-1}^{-2}.
 \end{align*}

Proposition~\ref{HICrit} now yields the theorem.

\begin{thm}
The space $\XK$ is hereditarily indecomposable.
\end{thm}

\section{Bounded linear operators on $\XK$}\label{Operators}

For technical reasons it will be convenient in the first few results of this
section to work with elements of $\XK$ all of whose coordinates are rational,
that is to say with elements of $\XK\cap \rats^\Gamma$.  Since (as may be
readily checked) each $d_\xi$ is in $\XK\cap \rats^\Gamma$, as are all rational
linear combinations of these, we see that $\XK\cap \rats^\Gamma$ is dense in
$\XK$.

\begin{lem}\label{RatVec}
Let $m<n$ be natural numbers and let $x\in \XK\cap \rats^\Gamma$, $y\in \XK$ be
such that $\ran x, \ran y$ are both contained in the interval $ (m,n]$. Suppose
that $\dist(y,\reals x)>\delta$.  Then there exists $b^*\in
\ball\ell_1(\Gamma_n\setminus \Gamma_m)$, with rational coordinates, such that
$b^*(x)=0$ and $b^*(y)>\frac12\delta$.
\end{lem}
\begin{proof}
Let $u,v\in \ell_\infty(\Gamma_n\setminus \Gamma_m)$ be the restrictions of
$x,y$ respectively.  Then $x=i_nu$, $y=i_nv$ and so, for any scalar $\lambda$,
$\|y-\lambda x\|\le \|i_n\|\|v-\lambda u\|$.  Hence $\dist(v,\reals
u)>\frac12\delta$ and so, by the Hahn--Banach Theorem in the finite dimensional
space $\ell_\infty(\Gamma_n\setminus \Gamma_m)$, there exists $a^*\in
\ball\ell_1(\Gamma_n\setminus \Gamma_m)$ with $a^*(u)=0$ and
$a^*(v)>\frac12\delta$. Since $x$ has rational coordinates our vector $u$ is in
$ \rats^{\Gamma_n\setminus \Gamma_m}$.  It follows that we can approximate
$a^*$ arbitrarily well with $b^*\in \rats^{\Gamma_n\setminus \Gamma_m}$
retaining the condition $b^*(u)=0$.
\end{proof}

\begin{lem}\label{DistExact}
Let $T$ be a bounded linear operator on $\XK$, let $(x_i)$ be a
$C$-RIS in $\XK\cap \rats^\Gamma$ and assume that
$\mathrm{dist}(Tx_i,\reals x_i)>\delta>0$ for all $i$.  Then, for
all $j,p\in \nats$, there exist $z\in [x_i:i\in \nats]$, $q>p$ and
$\eta\in \Delta_{q}$ such that
 \begin{enumerate}
 \item $(z,\eta)$ is a $(12C,2j,0)$-exact pair;
 \item $(Tz)(\eta)>\frac7{16}\delta$;
 \item $\|(I-P_{(p,q)})Tz\|<m_{2j}^{-1}\delta$;
 \item $\langle P^*_{(p,q]}e^*_{\eta},Tz\rangle >\frac{3}{8}\delta$.
 \end{enumerate}
\end{lem}
\begin{proof}
Since the sequence $(Tx_i)$ is weakly null, we may, by taking a subsequence if
necessary, assume that there exist $p<q_0<q_1<q_2<\dots$ such that, for all
$i\ge 1$, $\ran x_i\subseteq (q_{i-1},q_i)$ and
$\|(I-P_{(q_{i-1},q_i)})Tx_i\|<\frac15m_{2j}^{-2}\delta\le\frac1{80}m_{2j}^{-1}\le
\frac1{1280}\delta.$ It certainly follows from this that
$\mathrm{dist}(P_{(q_{i-1},q_i)}Tx_i,\reals x_i)>\frac{1279}{1280}\delta$. We
may  apply Lemma~\ref{RatVec} to obtain
$b^*_i\in\ball\ell_1(\Gamma_{q_i-1}\setminus\Gamma_{q_{i-1}})$, with rational
coordinates, satisfying
$$
\langle b^*_i,x_i\rangle=0,\quad \langle b^*_i,
P_{(q_{i-1},q_i)}Tx_i\rangle>\txtfrac{1279}{2560}\delta.
$$
Taking a further subsequence if necessary, we may assume that the coordinates
of $b^*_i$ have denominators dividing $N_{q_i-1}!$, so that $b^*_i\in
B_{q_i-1,q_{i-1}}$, and we may also assume that $q_1\ge 2j$. We are thus in a
position to apply Lemma~\ref{RISZeroExact}, getting elements $\xi_i$ of weight
$m_{2j}$ in $\Delta_{q_j}$ such that the element $\eta=\xi_{{n_{2j}}}$ of
$\Delta_{q_{n_{2j}}}$ has evaluation analysis
$$
e^*_\eta =
\sum_{i=1}^{n_{2j}}d^*_{\xi_i}+m_{2j}^{-1}\sum_{i=1}^{n_{2j}}P^*_{(q_{i-1},q_i)}b^*_i.
$$
and such that $(x,\eta)$ is a $(12C,2j,0)$-exact pair, where $z$
denotes the weighted average
$$
x=m_{2j}n_{2j}^{-1}\sum_{i=1}^{n_{2j}}x_i.
$$

We next need to estimate $(Tz)(\eta)$. For each $k$, we have
$\|(I-P_{(q_{k-1},q_k)})Tx_k\|<\frac1{80}m_{2j}^{-1}\delta$ so that
\begin{align*}
(Tx_k)(\eta) &\ge \langle e^*_\eta,
P_{(q_{k-1},q_k)}Tx_k\rangle-\txtfrac1{80}m_{2j}^{-1}\delta\\
&=m_{2j}^{-1}\langle b^*_k,
P_{(l_{k-1},l_k)}Tx_k\rangle-\txtfrac1{80}m_{2j}^{-1}\delta>\txtfrac{1247}{2560}m_{2j}^{-1}\delta.
\end{align*}
It follows that
 $$ (Tz)(\eta) =
 n_{2j}^{-1}m_{2j}\sum_{k=1}^{n_{2j}}(Tx_k)(\eta)>\frac{7}{16}\delta.
 $$

For inequality (3) in which we are taking $q=q_{n_{2j}}$, we note that
$p<q_{k-1}<q_k\le q$ for all $k$ so that
\begin{align*}
\|(I-P_{(p,q]})Tx_k\| &= \|(P_{(0,p]}+P_{(q,\infty)})Tx_k\|\\
&= \|(P_{(0,p]}+P_{(q,\infty)})(I-P_{(q_{k-1},q_k)})Tx_k\|\\
&\le 5\|(I-P_{(q_{k-1},q_k)})Tx_k\|< m_{2j}^{-2}\delta,
\end{align*}
using our usual estimates for norms of FDD projections.  The
inequality for the weighted average $z$  follows at once. Inequality
(4) follows  from (2) and (3) thus
$$
\langle P^*_{(p,q]}e^*_{\eta},Tz\rangle \ge (Tz)(\eta) =
\|(I-P_{(p,q]})Tz\|
> \txtfrac7{16}\delta-m_{2j}^{-1}\delta\ge \txtfrac3{16}\delta.
$$
\end{proof}

\begin{prop}\label{distLem}
Let $T$ be a bounded linear operator on $\XK$ and let $(x_i)_{i\in \nats}$ be a
RIS in $\XK$.  Then $\mathrm{dist}(Tx_i,\reals x_i)\to 0$ as $i\to \infty$.
\end{prop}
\begin{proof}
It will be enough to prove the result for a RIS in $\XK\cap \rats^\Gamma$.
Suppose, if possible, that $\mathrm{dist}(Tx_i,\reals x_i)>\delta>0$ for all
$i$.  The idea is to  obtain a dependent sequence in rather the same way as we
did in Lemma~\ref{HILem}, except that this time it will be a 0-dependent
sequence, rather than a 1-dependent sequence.

We start by choosing $j_0$ such that $m_{2j_0-1}> 256C\|T\|\delta^{-1}$ and
$j_1$ such that $m_{4j_1-2}>n_{2j_0-1}^2$. Taking $p=p_0=0$ and $j=2j_1-1$ in
Lemma~\ref{DistExact} we can find $q_1$ and a $(12C,4j_1-2,0)$-exact pair
$(z_1,\eta_1)$ with $\rank \eta_1=q_1$, $(Tz_1)(\eta_1)>\frac38\delta$ and
$\|(I-P_{(0,q_1]})(Tz_1)\|<m_{4j_1-2}^{-1}\delta$. Let $p_1=q_1+1$ and let
$\xi_1$ be the special Type  1 element of $\Delta_{p_1}$ given by  $\xi_1 =
(p_1,m_{2j_0-1},e^*_{\eta_1})$.

Now, recursively for $2\le i\le n_{2j_0-1}$, define $j_{i}=\sigma(\xi_{i-1})$,
 and use the lemma again to choose $q_{i}$ and a
$(12C,4j_{i},0)$-exact pair $(z_{i},\eta_{i})$ with $\rank \eta_i=q_i$, $\ran
z_i\subseteq (p_{i-1},q_i]$, $\langle
P^*_{(p_{i-1},q_i]}e^*_{\eta_i},Tz_1\rangle>\frac38\delta$ and
$\|(I-P_{(p_i,q_i]})(Tz_i)\|<m_{4j_i}^{-1}\delta$.
 We now define $p_i=q_i+1$ and let $\xi_i$ to be the Type 2 element
$(p_i,\xi_{i-1},m_{2j_0-1}^{-1}, e^*_{\eta_i})$ of $\Delta_{p_i}$.

It is clear that we have constructed a $(12C,2j_0-1,0)$-dependent sequence
$(z_i)_{1\le i\le n_{2j_0-1}}$.

By the estimate of Lemma~\ref{ZeroDepSeqEst} we have
$$
\|z\|\le 48Cm_{2j_0-1}^{-2}
$$
 for the average
$$
z=n_{2j_0-1}^{-1}\sum_{i=1}^{n_{2j_0-1}}z_i.
$$
However, let us consider the element $\gamma=\xi_{n_{2j_0-1}}$ of
$\Delta_{p_{n_{2j_0-1}}}$, which has evaluation analysis
$$
e^*_\gamma = \sum_{i=1}^{n_{2j_0-1}} d^*_{\xi_i} +
m_{2j_0-1}^{-1}\sum_{i=1}^{n_{2j_0-1}}P^*_{(p_{i-1},p_i)}e^*_{\eta_i}.
$$
Noting that $p_k=q_q+1$ for $k\ge 1$, and that
$m_{4j_i}>m_{4j_1-2}>n_{2j_0-1}^2$, we may estimate $(Tz)(\gamma)$ as follows
\begin{align*}
(Tz)(\gamma) &= n_{2j_0-1}^{-1}\sum_{k=1}^{n_{2j_0-1}}(Tz_k)(\gamma)\\
&\ge   n_{2j_0-1}^{-1}\sum_{k=1}^{n_{2j_0-1}}\left(\langle P^*_{(p_{k-1},p_{k})}e^*_\gamma,Tz_k\rangle-\|(I-P_{(p_{k-1},q_k]]})(Tx_k)\|\right)\\
&\ge n_{2j_0-1}^{-1} \sum_{k=1}^{n_{2j_0-1}} \left(m_{2j_0-1}^{-1}\langle P^*_{(p_{k-1},p_{k})}e^*_{\eta_k},Tx_k\rangle- m_{4j_1-2}^{-1}\delta\right)\\
&\ge \delta n_{2j_0-1}^{-1}
\sum_{k=1}^{n_{2j_0-1}}(\txtfrac38m_{2j_0-1}^{-1}-5n_{2j_0-1}^{-2})>\txtfrac14m_{2j_0-1}^{-1}\delta.
 \end{align*}
 So
 $$
 \|Tz\| \ge  \txtfrac14m_{2j_0-1}^{-1}>\txtfrac1{144}C^{-1}\delta m_{2j_0-1}\|z\|,
 $$
which  is a contradiction because $\txtfrac1{144}C^{-1}\delta m_{2j_0-1}>\|T\|$
by our original choice of $j_0$.
 \end{proof}

\begin{thm} Let $T$ be a bounded linear operator on $\XK$.  Then
there exists a scalar $\lambda$ such that $T-\lambda I$ is compact.
\end{thm}
\begin{proof}
We start by considering a normalized RIS $(x_i)$ in $\XK$.  By
Proposition~\ref{distLem} there exist scalars $\lambda_i$ such that
$\|Tx_i-\lambda_ix_i\|\to 0$.  We claim that $\lambda_i$ necessarily tends to
some limit $\lambda$.  Indeed, if not, by passing to a subsequence,  we may
suppose that $|\lambda_{i+1}-\lambda_{i}|>\delta$ for all $i$. Now the sequence
$(y_i)$ where $y_i= x_{2i-1}+x_{2i}$ is again a RIS, so that there exist
$\mu_i$ with $\|Ty_i-\mu_iy_i\|\to 0$ by Proposition~\ref{distLem} again. We
thus have
$$
\|(\lambda_{2i}-\mu_i)x_{2i} + (\lambda_{2i-1}-\mu_i)x_{2i-1}\|\le
\|Tx_{2i}-\lambda_{2i}x_{2i}\|+\|Tx_{2i-1}-\lambda_{2i-1}x_{2i-1}\|+
\|Ty_i-\mu_iy_i\|\to 0.
$$
Since the RIS $(x_i)$ is a block sequence, there exist $l_i$ such that
$P_{(0,l_i]}y_i = x_{2i-1}$ and $P_{(l_i,\infty)}y_i = x_{2i}$.  Using the
assumption that the sequence $(x_i)$ is normalized we now have
$$
|\lambda_{2i-1}-\mu_i|= \|(\lambda_{2i-1}-\mu_i)x_{2i-1}\| \le
\|P_{(0,l_i]}\|\|(\lambda_{2i}-\mu_i)x_{2i} + (\lambda_{2i-1}-\mu_i)x_{2i-1}\|,
$$
with a similar estimate for $|\lambda_{2i}-\mu_i|$.  Each of these sequences
thus tends to 0, so that $\lambda_{2i}-\lambda_{2i-1}$ also tends to 0,
contrary to our assumption.

We now show that the scalar $\lambda$ is the same for all rapidly increasing
sequences.  Indeed, if $(x_i)$ and $(x'_i)$ are RIS with $\|Tx_i-\lambda
x_i\|\to 0$ and $\|Tx'_i-\lambda' x_i\|\to 0$, we may find $i_1<i_2<\cdots$
such that the sequence $(y_k)$ defined by
$$
y_k = \begin{cases} x_{i_k}\ \quad\text{if $k$ is odd}\\
                    x'_{i_k}\ \quad\text{if $k$ is even}
    \end{cases}
$$
is again a RIS.  By the first part of the proof we must have
$\lambda=\lambda'$.

We have now obtained $\lambda$ such that $\|(T-\lambda I)x_i\|\to 0$ for every
RIS.  By Proposition~\ref{RIStoBlock}, we deduce that $\|(T-\lambda I)x_i\|\to
0$ for every bounded block sequence in $\XK$.  This, of course, implies that
$T-\lambda I$ is compact.
\end{proof}

\section{Reflexive subspaces with HI duals}\label{HIDuals}

We devote this section to a proof that $\XK$ is saturated with reflexive HI
subspaces having HI duals.  The proof involves reworking much of the
construction of Section~\ref{ellOneExact} in the context of a subspace of $\XK$
and its dual. By standard blocking arguments, it is enough to prove the
following theorem.

\begin{thm}\label{HIDualThm}
Let $L=\{l_0,l_1,l_2,\dots\}$ be a set of natural numbers satisfying
$l_{n-1}+1<l_{n}$, and for each $n\ge 1$ let $F_n$ be a subspace of the
finite-dimensional space $P_{(l_{n-1},l_{n})}\XK=\bigoplus_{l_{n-1}<k<l_{n}}
M_k$.  Then the subspace $W= \overline{\bigoplus_{n\in\nats} F_n}$ of $\XK$ is
reflexive and has HI dual.
\end{thm}

We note in passing the following corollary, which gives an
indication of the ``very conditional'' nature of the basis of
$\ell_1$ that we have constructed.  For the purposes of the
statement we briefly abandon the ``$\Gamma$ notation'' and revert to
the notation of Definition~\ref{TriangDef} and Theorem~\ref{BDThm}.

\begin{cor}
There exist a basis $(d^*_n)_{n\in \mathbb N}$ of $\ell_1$ and
natural numbers $k_1<k_2<\cdots$ with the property that the quotient
$\ell_1/[d^*_n:n\in M]$ is hereditarily indecomposable whenever the
subset $M$ of $\mathbb N$ contains infinitely many of the intervals
$(k_p,k_{p+1}]$.
\end{cor}

The rest of this section will be devoted to the proof of
Theorem~\ref{HIDualThm}. We have already remarked at the end of
Section~5 that the subspace $W$ defined in the statement of the
theorem is reflexive. The subspaces $F_n$ form a finite-dimensional
decomposition of $W$, the corresponding FDD projections being
$Q_{(m,n]}=P_{(l_m,l_n]}\restriction W= P_{(l_m,l_n)}\restriction
W$, when $0\le m<n$.  The dual space $W^*$ has a dual FDD $(F_n^*)$
and corresponding projections $Q^*_{(m,n]}$. We shall establish
hereditary indecomposability of $W^*$ via the criterion
Proposition~\ref{HICrit}.  We write $R$ for the quotient mapping
$\XK^*=\ell_1\to W^*$ and observe that if $f^*_n\in F_n^*$ for $1\le
n\le N$ then the norm of $f^*=\sum_{n=1}^N f^*_n$ in $W^*$ is given
by
$$
\|\sum f^*\|_{Y^*} = \inf\{\|g^*\|:g^*\in \XK^*\text{ and } Rg^*=f^*\}.
$$

\begin{lem}\label{QuotSupp}
If $f^*\in \im Q^*_{(M,N]}=\bigoplus_{M<n\le N}F^*_n\subset W^*$ then there
exists $h^*\in \XK^*=\ell_1(\Gamma)$ with $\supp h^*\subseteq
\Gamma_{l_N-1}\setminus \Gamma_{l_M}$ and $\|h^*\|_1\le 4\|f^*\|$ and
$RP^*_{(l_M,l_N)}h^*=RP^*_{(l_M,\infty)}h^*=f^*$.
\end{lem}
\begin{proof}
We extend $f^*$ by the Hahn--Banach theorem to obtain $g^*\in
\XK^*=\ell_1^(\Gamma)$ with $Rg^*=f^*$ and $\|g^*\|_{\XK^*}=\|f^*\|_{W^*}$. We
set $h^*=P_{(0,l_N)}g^*\in \ell_1(\Gamma_{l_N-1})$ and
$b^*=h^*\chi_{\Gamma_{l_N-1}\setminus \Gamma_{l_N}}$, noting that
$$
\|b^*\|_1\le \|h^*\|_1\le 2\|g^*\|_1\le 4\|g^*\|_{\XK^*}=4\|f^*\|.
$$
To check that $RP^*_{(l_M,l_N)}h^*=RP^*_{(l_M,\infty)}b^*=f^*$, we first note
that
$$
P^*_{(l_M,\infty)}b^* =  P^*_{(l_M,\infty)}h^*,
 $$
because $P^*_{(l_M,\infty)}k^*=0$ whenever $\supp k^*\subseteq
\Gamma_{l_M}$. Since both $b^*$ and $h^*$ are supported by
$\Gamma_{l_N-1}$ we have
$$
P^*_{(l_M,l_N)}b^* =
P^*_{(l_M,\infty)}P^*_{(0,l_N)}b^*=P^*_{(l_M,\infty)}b^*=P^*_{(l_M,\infty)}h^*=P^*_{(l_M,\infty)}P^*_{(0,l_N)}h^*
= P^*_{(l_M,l_N)}g^*.
 $$
It follows that
$$
R^*P^*_{(l_M,l_N)}b^*=R^*P^*_{(l_M,l_N)}g^*=g^*\circ P_{(l_M,l_N)}\restriction
W= g^*\circ Q_{(M,N]}=f^*.
$$
\end{proof}

\begin{lem}\label{DualUpperEst} Let $j\ge 1$, $1\le a\le n_{2j}$ and  $M\le M_0<M_1<M_2<\cdots<M_{a}$ be natural
numbers, with $2j\le M_1$. For each $i\le a$, let $f^*_i$ be in $\ball
\bigoplus_{M_{i-1}<n\le M_i} F^*_n$ and write $f^*=\sum_{i=1}^{a}f^*_i$. Then
there exists $\gamma\in \Gamma$ with $p^*_{(0,l_M]}e^*_\gamma=0$ and
$\|4m_{2j}R(e^*_\gamma)-f^*\|\le 2^{-l_M+3}$; in particular $\|f^*\|_{Y^*}\le
5m_{2j}$.
\end{lem}
\begin{proof}
By Lemma~\ref{QuotSupp} there exist $h^*_i\in
\ell_1(\Gamma_{l_{M_i-1}}\setminus \Gamma_{l_{M_{i-1}}})$ with $\|h^*_i\|_1\le
4$ and $R(P^*_{(l_{M_{i-1}},l_{M_i}i})=f^*_i$. Since
$B_{l_{M_i}-1,l_{M_{i-1}}}$ is an $\epsilon$-net in
$\ball\ell_1(\Gamma_{l_{M_i}-1}\setminus \Gamma_{l_{M_{i-1}}})$, with $\epsilon
= 2^{-l_{M_i}+1}\le 2^{-l_M-2i+1}$ we can choose $b^*_i\in
B_{l_{M_k},l_{M_{k-1}}}$ such $\|h^*_i-4b^*_i\|_1\le 2^{-l_M-2i+3}$.

Now write $p_i=l_{M_i}$ for $1\le i\le a$ and apply the construction of
Proposition~\ref{ExistGammaEven} to obtain $\gamma\in \Delta_{p_a}$ with
evaluation analysis
$$
e^*_\gamma = \sum_{i=1}^{a} d^*_{\xi_i} + m_{2j}^{-1}\sum_{i=1}^{a}
P^*_{(p_{k-1},\infty)}b^*_k.
$$
Since $\rank \xi_i =p_i\in  L$ for all $i$, we have $Rd^*_{\xi_i}=0$ and so
\begin{align*}
\|f^*-2m_{2j}R(e^*_\gamma)\| &= \|\sum_{i=1}^{a} \bigl
(f^*_i-2RP^*_{(p_{i-1},\infty)}b^*_i\bigr)\|\\
&\le \sum_{i=1}^{n_{2j}}
\|RP^*_{(p{i-1},\infty]}h^*_i-2RP^*_{(p_{i-1},\infty)}b^*_i\|\\
 &\le 3\sum_{i=1}^{a}
\|h^*_i-2b^*_i\| \le 3\sum_{i=1}^{\infty} 2^{-l_M-2i+2}=2^{-l_M+2}.
\end{align*}
It follows that $\|f^*\|\le \|4m_{2j}R(e^*_\gamma)\|+8\le 5m_{2j}$.
\end{proof}

\begin{lem} \label{SubspAve} Let $Y$ be any block subspace  of $W^*$
and let $n,M$ be positive integers. For every $C>1$ there exists a
$4C$-$\ell_1^n$-average $w\in W$, with $Q_{(0,M]}w=0$, and a functional $g^*\in
\ball Y$ with $Q^*_{(0,M]}g^*=0$  and $\langle g^*,w\rangle \ge 1$.
\end{lem}
\begin{proof}
The proof is a dualized version of Lemma~\ref{ExistAve}.  We suppose, without
loss of generality, that $C<2$ and choose $l,j$ such that $C^l>n$ and
$n_{2j}>(10n_{2j})^l$; we take $k$ minimal subject to $C^k>5m_{2j}$ noting that
$$
n_{2j}>(10m_{2j})^l\ge (2C^{k-1})^l\ge C^{kl}>n^k.
$$
Now take $(f^*_i)_{i=1}^{n^k}$ to be a normalized block sequence in $Y\cap \ker
Q^*_{(0,M]}$; we may apply Lemma~\ref{DualUpperEst} to obtain
$$
\|\sum_{i=1}^{n^k} \pm f^*_i\| \le 5m_{2j}<C^k.
$$
So by part (ii) of Lemma~\ref{ellOneAve} (with $C=1+\epsilon$) there are
successive linear combinations $g^*_1,\dots, g^*_{n}$ such that $\|g^*_i\|\ge
C^{-1}$ for all $i$, while
$$
\|\sum_{i=1}^n \pm g^*_i\|\le 1,
$$
for all choices of sign. Since $(g^*_i)$ is a block sequence in $\ker
Q^*_{(0,M]}$ we can choose $M\le N_0<N_1<\dots$ such that
$Q^*_{(N_{i-1},N_i]}g^*_i=g^*_i$.  Now we choose, for each $i$ an element $w_i$
of $W$ such that $\|w_i\|\le C$ and $\langle g^*_i,w_i\rangle =1$.  If we set
$w'_i =Q_{(N_{i-1},N_i]}w_i$ then we have $\|w'_i\|\le 4C$ and $\langle
g^*_i,w'_i\rangle =\langle g^*_i,w_i\rangle =1$, while $\langle
g^*_i,w'_h\rangle=0$ when $h\ne i$. The element $w=n^{-1}\sum_{i=1}^n w'_i$ is
thus a $4C$-$\ell_1^n$ average, with $Q_{(0,p]}w=0$, and satisfies $\langle
g^*,w\rangle = 1$, where $g^*= \sum_{i=1}g^*_i\in \ball Y$.
\end{proof}

\begin{lem}\label{DualExact}
Let $Y$ be any block subspace of $W^*$ and let $N,j$ be  positive integers.
There exists a $(600,2j,1)$-exact pair $(z, \gamma)$ with $z\in W$,
$Q_{(0,N]}z=0$, $P^*_{(0,l_N]}e^*_\gamma=0$ and
$\dist(Re^*_\gamma,Y)<2^{-l_N}$.
\end{lem}
\begin{proof}
By repeated applications of Lemma~\ref{SubspAve}, we construct natural numbers
$N\le M_0<M_1<M_2<\cdots$ and $j_1<j_2<\dots$, elements
$w_i=Q_{(M_{i-1},M_i]}w_i$ of $W$, and functionals
$g^*_i=Q^*_{(M_{i-1},M_i]}g^*_i\in \ball Y$ such that
\begin{enumerate}
\item $w_i$ is a $5$-$\ell_1^{n_{j_i}}$-average;
\item $\langle g^*_i,w_i\rangle \ge 1$;
\item $j_{i+1}>M_i$.
\end{enumerate}
It follows from Lemma~\ref{AveRIS} that $(w_i)$ is a $10$-RIS.

Writing $g^*=\sum_{i=1}^{n_{2j}}g^*_i$ and applying Lemma~\ref{DualUpperEst} we
find $\gamma$ of weight $m_{2j}$ such that $\|4m_{2j}R(e^*_\gamma)-g^*\|\le
2^{-N+3}$.  We thus have
$$
\dist(Re^*_\gamma,Y)\le \|Re^*_\gamma-\txtfrac14m_{2j}^{-1}g^*\|\le
2^{-l_N+1}m_{2j}^{-1}<2^{-l_N},
$$
and
$$
4m_{2j}\sum_{i=1}^{n_{2j}}w_i(\gamma) \ge \sum_{i=1}^{n_{2j}}\langle
g^*,w_i\rangle -2^{-l_N+3}\ge n_{2j}-16.
$$
We now set $z= \theta m_{2j}n_{2j}^{-1}\sum_{i=1}^{n_{2j}}w_i $ where $\theta$
is chosen so that $z(\gamma)=1$; by the above inequality $0<\theta\le
4+128n_{2j}^{-1}<5$.

To estimate $\|z\|$ and $|z(\gamma')|$ when $\weight \gamma'= m_h\ne m_{2j}$ we
return to Lemma~\ref{RISUpperEst} deducing that
$$
\|z\|\le 60\theta\quad\text{and}\quad|z(\gamma')| \le \begin{cases} 110\theta
m_{h}^{-1} \quad\text{ if }
h<2j\\
               \ \, 60\theta m_{2j}^{-1} \quad\text{ if } h>2j.
                \end{cases}
$$
So $(z, \gamma)$ is certainly a $(600,2j,1)$-exact pair.
\end{proof}

\begin{lem}\label{DualHILem}
Let $Y_1$ and $Y_2$ be block subspaces of $W^*$ and let $j_0$ be a natural
number. There exists a sequence $(x_i)_{i\le n_{2j_0-1}}$ in $W$, together with
natural numbers $0=p_0<p_1<p_2<\cdots<p_{n_{2j_0-1}}$, and elements
$\eta_i\in\Gamma_{p_i-1}\setminus \Gamma_{p_{i-1}}$, $\xi_i\in \Delta_{p_i}$
($1\le i\le n_{2j_0-1}$), satisfying the conditions {\rm (1)} to {\rm (4)} of
Definition~\ref{DepSeq} with $C=600$, $\varepsilon =1$, and such that, for all
$i\ge 1$, the following additional properties hold
\begin{enumerate}
\setcounter{enumi}4 \item $\rank \xi_i = p_i\in L$;
\item $P^*_{(p_{i-1},p_i]}e^*_{\eta_i}=0$, $P_{(p_{i-1},p_i]}(x_i)=x_i$;
\item $\dist (Re^*_{\eta_i},Y_k) < 2^{-p_{i-1}}$,
where $k=1$ for odd $i$ and $k=2$ for even $i$.
\end{enumerate}
\end{lem}
\begin{proof}
We start by choosing $j_1$ such that $m_{4j_1-2}>n_{2j_0-1}^2$ and then
applying Lemma~\ref{DualExact} to obtain a $(600,4j_1-2,1)$-exact pair
$(x_1,\eta_1)$ with $x_1\in W$. Set $p_1=l_{N_1}$, where $N_1$ is large enough
to ensure that $P_{(0,p_1)}x_1=Q_{(0,N_1]}x_1=x_1$, $\rank \eta_1<p_1$ and
$2^{p_1}>2n_{2j_0-1}$. Let $\xi_1=(p_1,m_{2j_0-1}^{-1},\eta_1)\in
\Delta_{p_1}$.

Continuing recursively, if for some $i<n_{2j_0-1}$, we have defined $\xi_i\in
\Delta_{p_i}$, where $p_i=l_{N_i}$, we set $j_{i+1}=\sigma(\xi_i)$ and apply
Lemma~\ref{DualExact} to get a $(600,4j_{i+1},1)$-exact pair
$(x_{i+1},\eta_{i+1})$ with $x_{i+1}\in W$,
$Q_{(0,N_{i}]}x_{i+1}=P_{(0,p_i]}x_{i+1}=0$, $P^*_{(0,p_i]}e^*_{\eta_{i+1}}=0$
and $\dist(R^*e^*_{\eta_{i+1}},Y_k)<2^{-p_i}$, where $k$ depends on the parity
of $i+1$. We now take $N_{i+1}$ large enough, set $p_{i+1}=l_{N_{i+1}}$ and
define $\xi_{i+1}=(p_{i+1},\xi_i,m_{2j_0-1}^{-1}, \eta_{i+1})\in
\Delta_{p_{i+1}}$.
\end{proof}

We are now ready to finish the proof of the theorem.  We consider any two
infinite-dimensional subspaces $Y_1$ and $Y_2$ of $W^*$ and apply
Lemma~\ref{DualHILem} obtaining a dependent sequence satisfying (1) to (7).  By
property (7) we may choose, for each $i$, an element $y^*_i$ of $Y_k$ with
$$
\|y^*_i-Re^*_{\eta_i}\|<2^{-p_{i}}.
$$
We set
$$
y^*=m_{2j_0-1}^{-1}\sum_{i\text{ odd}} y^*_i\in Y_1,\quad
z^*=m_{2j_0-1}^{-1}\sum_{i\text{ even}} y^*_i\in Y_2.
$$
If $\gamma$ is the element $\xi_{n_{2j_0-1}}$ then the evaluation analysis of
$\gamma$ is
\begin{align*}
e^*_\gamma &= \sum_{i=1}^{n_{2j_0-1}}d^*_{\xi_i} +
m_{2j_0-1}^{-1}\sum_{i=1}^{n_{2j_0-1}}p^*_{(p_{i-1},\infty)}e^*_{\eta_i}\\
&=\sum_{i=1}^{n_{2j_0-1}}d^*_{\xi_i} +
m_{2j_0-1}^{-1}\sum_{i=1}^{n_{2j_0-1}}e^*_{\eta_i},
\end{align*}
because $P^*_{(0,p_{i-1}]}e^*_{\eta_i}=0$.  Since $\rank \xi_i=p_i\in L$ for
all $i$ we have
$$
Re^*_\gamma = m_{2j_0-1}^{-1}\sum_{i=1}^{n_{2j_0-1}}Re^*_{\eta_i},
$$
which leads to
$$
\|y^*+z^*\|\le 1+ \|m_{2j_0-1}^{-1}\sum_{i=1}^{n_{2j_0-1}}Re^*_{\eta_i}\|=
1+\|Re^*_\gamma\|\le2.
$$

We shall prove that $\|y^*-z^*\|$ is very large by estimating
$\langle y^*-z^*,x\rangle$, where $x$ is the average
 $$
  x =n_{2j_0-1}^{-1}\sum_{k=1}^{n_{2j_0-1}}(-1)^kx_k,
  $$
about which we know from Lemma~\ref{OneDepSeqEst} that
$$
\|x\|\le 7200m_{2j_0-1}^{-2}.
$$
By (7) and the definition of a 1-exact pair, we have
 $$
 \langle e^*{\eta_i},x_k\rangle =\begin{cases} 1\quad\text{if }i=k\\
  0\quad\text{if }i\ne k,\end{cases}
   $$ so that
\begin{align*}
\langle y^*-z^*,x\rangle &= n_{2j_0-1}^{-1}m_{2j_0-1}^{-1} \sum{i,k}\langle y^*_i,x_k\rangle\\
&\ge n_{2j_0-1}^{-1}m_{2j_0-1}^{-1}\sum_{i,k}(\langle
e^*{\eta_i},x_k\rangle-2^{-p_i})\\
&\ge m_{2j_0-1}^{-1}(1 - n_{2j_0-1}2^{-p_1})\ge \txtfrac12m_{2j_0-1}^{-1},
\end{align*}
the last step following from our choice of $p_1$ with $2^{p_1}>2n_{2j_0-1}$.

We can now deduce that
$$
\|y^*-z^*\|\ge \frac{m_{2j_0-1}}{14400}.
$$
We have shown that the subspaces $Y_1$ and $Y_2$ of $W^*$ contain
elements $y^*$, $z^*$ with $\|y^*+z^*\|\le 2$ and $\|y^*-z^*\|$
arbitrarily large.  By Proposition~\ref{HICrit}, we have established
hereditary indecomposability of $W^*$.

\section{Concluding Remarks}
\label{ConcRemSection}

\subsection{Operators on subspaces of $\XK$}

If we are looking at a bounded linear operator $T:Y\to \XK$ defined only on a
subspace $Y$ of $\XK$, rather than on the whole space, then, as in other HI
constructions, the arguments of the preceding section can be used to show that
$T$ can be expressed as $\lambda I_Y+S$ with $S$ strictly singular.  However,
as we shall now see, in this case the perturbation need not be compact.

\begin{prop}
There exists a subspace $Y$ of $\XK$ and a strictly singular, non-compact
operator $T$ from $Y$ into $\XK$.  In fact, for a suitably chosen $Y$, we may
choose $T$ mapping $Y$ into itself.
\end{prop}
\begin{proof}
By a theorem of Androulakis, Odell, Schlumprecht and Tomczak-Jaegermann
\cite{AOST}, in order to find $Y$ and a strictly singular, non-compact $T:Y\to
\XK$,  it is enough to exhibit normalized sequences $(x_i)$ and $(y_i)$ in
$\XK$ such that $(y_i)$ has a spreading model equivalent to the usual
$\ell_1$-basis, while $(x_i)$ has a spreading model that is not equivalent to
that basis. For $(x_i)$ we may take any normalized RIS; indeed, by
Proposition~\ref{BasicIneq}, the spreading model associated with any RIS  is
dominated by the unit vector basis of the Mixed Tsirelson space $\mathfrak
T[(\mathscr A_{3n_j},m_j^{-1})_{j\in \nats}]$, and so is not equivalent to the
$\ell_1$-basis.  For $(y_i)$ we may take a specific sequence , setting
$$
y_n = \sum_{\xi\in \Delta_n} d_\xi.
$$
The result we need is a lemma about norms of linear combinations of these
vectors.

\begin{lem}
Let $F$ be a finite set of natural numbers with $\min F\ge j$ and $\#F<
2n_{2j}$.  Then, for all real scalars $a_n$,
$$
\|\sum_{n\in F} a_n y_n\|\ge \txtfrac14\sum_{n\in F}|a_n|.
$$
\end{lem}
\begin{proof}
Without loss of generality, we may suppose that $\sum_{n\in F} a^+_n\ge
\txtfrac12\sum_{n\in F}|a_n|$ and we may choose $p_1,p_2,\dots,p_r$ in $F$,
with $p_{i+1}>p_i+1$, $r\le p_{2j}$, and
$$
\sum_{i=1}^r a_i \ge \txtfrac14\sum_{n\in F}|a_n|.
$$

Since $p_1\ge \min F\ge 2j$, $\Delta_{n_1}$ does contain Type 1 elements of the
form $(p_1,m_{2j}^{-1},\pm e^*_{\eta_1})$, with $\eta_1\in \Gamma_{n_1-1}$. We
take $\xi_1$ to be such an element, and continue recursively, for $1\le i<r$,
taking $\eta_{i+1}$ to be any element of $\Delta_{p_i+1}$ and $\xi_{i+1}$ to be
the Type 2 element $(p_{i+1},\xi_i,m_{2j}^{-1},\pm e^*_{\eta_{i+1}})$ of
$\Delta_{n_{i+1}}$. If $\gamma=\xi_{r}$ then the evaluation analysis of
$\Gamma$ is
$$
e^*_\gamma = \sum_{i=1}^r d^*_{\xi_i} + m_{2j}^{-1}\sum_{i=1}^r \pm
P^*_{(n_{i-1},n_i)}e^*_{\eta_i}.
$$
If we write $y=\sum_{n\in F} a_ny_n$, we have $\langle
d^*_{\xi_i},y\rangle=a_{n_i}$ for each $i$, so that
$$
e^*_\gamma(y) =\sum_{i=1}^r a_{p_i} + m_{2j}^{-1}\sum_{i=1}^r \pm
P^*_{(n_{i-1},n_i)}e^*_{\eta_i}(y).
$$
We have not until now been explicit about how the signs $\pm$ were chosen, but
it is now clear that this may be done in such a way that $e^*_\gamma(y)\ge
\sum_{i=1}^r a_{p_i}\ge \txtfrac14\sum_{n\in F}|a_n|.$
\end{proof}

It is now clear that the theorem of Androulakis et al may be applied.  In order
to get the refined version where $T$ takes $Y$ into  itself, it is enough to
look a little more closely at the proof given in \cite{AOST}. It turns out that
we may take $(y_i)$ as above and $Y$ to be the closed linear span $[y_i:i\in
\nats]$. It may be shown that, for any RIS $(x_i)$, the mapping $y_i\mapsto
x_i$ extends to a bounded linear operator from $Y$ to $\XK$. Since $Y$, like
all other infinite dimensional subspaces, contains a RIS, we may choose the
$x_i$ to lie in $Y$.
\end{proof}

\subsection{Very incomparable Banach spaces}

The original spaces $X_{a,b}$ of Bourgain and Delbaen provided, for the first
time, a continuum of non-isomorphic $\mathscr L_\infty spaces$.  It has also
been noted \cite{AAF} that if we take $Y$ to be Hilbert space and $X$ to be
$X_{a,b}$ with (for instance) $0<b<\frac12<a<1$, $a^4+b^4=1$, then all
operators from $X$ to $Y$ and all operators from $Y$ to $X$ are compact. The
constructions in the present paper allow us to exhibit a continuum of spaces
$X_\alpha$ ($\alpha\in \mathfrak c$) such that, for all $\alpha\ne \beta$, $
\mathcal L(X_\alpha,X_\beta)=\mathcal K(X_\alpha,X_\beta)$.

We start by taking an almost-disjoint family $(L_\alpha)_{\alpha\in \mathfrak
c}$ of  infinite subsets of $\nats$.  For each $\alpha$ we enumerate $L_\alpha$
in increasing order as $l^\alpha_j$ and define
$$
m^\alpha_j = m_{l^\alpha_j}, \qquad n^\alpha_j = n_{l^\alpha_j},
$$
where $(m_j,n_j)= (2^{2^j}, 2^{2^{j^2+1}})$ is the sequence mentioned in
Subsection~\ref{MT}.

Now  we may take $X_\alpha$ to be either  $\BmT[(\mathscr A_{n^\alpha_j},
1/m^\alpha_j)_{j\in \nats}]$ or $\XK[(\mathscr A_{n^\alpha_j},
1/m^\alpha_j)_{j\in \nats}]$.

\begin{prop}
Assume that $\alpha\ne \beta$ and let $T:X_\alpha\to X_\beta$ be  a bounded
linear operator.  For any RIS $(x_i)_{i\in \nats}$ in $X_\alpha$, we have
$\|T(x_i)\|\to 0$ as $i\to \infty$.
\end{prop}
\begin{proof}
Let $(x_i)$ be a $C$-RIS in $X_\alpha$ and suppose, if possible, that
$\|Tx_i\|>\delta>0$ for all $i$. Since $(Tx_i)$ is weakly null we may, by
taking a subsequence, assume that $(Tx_i)$ is a small perturbation of a
skipped-block sequence in $X_\beta$.  Thus, if $l=l^\beta_{2j}\in L_\beta$, we
may apply Proposition~\ref{LowerEst} to conclude
$$
\|n_{l}^{-1}\sum_{i=1}^{n_l} Tx_r\|_{X_\beta}\ge \txtfrac14 m_{2j}^{-1}n_l^{-1}
\sum_{r=1}^{n_l} \|Tx_r\|\ge \txtfrac14 \delta m_{2j}^{-1}.
 $$
On the other hand, Corollary~\ref{RISDom} tells us that
$$
\|n_{l}^{-1}\sum_{i=1}^{n_l} x_r\|_{X_\alpha}\le 10C
\|n_{l}^{-1}\sum_{i=1}^{n_l}e_i\|,
$$
where the norm on the right-hand side is calculated in $T[(\mathscr
A_{3n_j},m_j^{-1})_{j\in L_\alpha}]$.  If $l$ is not in $L_\alpha$ then this
norm is at most $m_l^{-2}$ by Lemma~\ref{MTLem}, so that
$$
\|n_{l}^{-1}\sum_{i=1}^{n_l} x_r\|_{X_\alpha}\le 10Cm_l^{-2}.
$$
By the assumed almost-disjointness of $L_\beta $ and $L_\alpha$ we can
certainly choose $j$ such that $l^\beta_{2j}\notin L_\alpha$ and $m_{l}>
40\|T\|\delta^{-1}$, yielding a contradiction.
\end{proof}

\begin{remark} The topologies $\sigma(\ell_1,X_\alpha)$ provide a continuum of
very incomparable weak$^*$ topologies on $\ell_1$: indeed, any linear mapping
on $\ell_1$ which is continuous from $\sigma(\ell_1,X_\alpha)$ to
$\sigma(\ell_1,X_\beta)$, with $\alpha \ne \beta$ is necessarily compact.
\end{remark}

\subsection{The space of operators $\mathcal L(\XK)$}

Of course, the spaces $\mathcal L(X)$ and $\mathcal K(X)$ of bounded
(respectively compact) linear operators on an infinite-dimensional
Banach space $X$ are always decomposable. (Indeed, for finite
dimensional subspaces $E\subset X$ and $F\subset X^*$, the subspaces
$X^*\otimes E$ and $F\otimes X$ are complemented.) So we must not
hope for too much exotic structure in these spaces of operators. In
this section we shall look briefly at subspaces of $\mathcal
L(\XK)$.  Certainly, $\mathcal L(\XK)=\mathcal K(\XK)\oplus \reals
I$ has HI subspaces, such as those isomorphic to $\XK$, and
subspaces isomorphic to $\XK^*=\ell_1$. It has no subspace
isomorphic to $c_0$ by a result of Emmanuele. (The main result of
\cite{E} shows that $c_0$ does not embed into $\mathcal K(X_{a,b})$
and the same proof works for $\XK$.) We shall now see that
$\mathcal(\XK)$  does have other subspaces with unconditional basis.
It is a general fact that if $(x_n)$ is a basic sequence in a Banach
space $X$ then the injective tensor product
$\ell_1\hat\otimes_\varepsilon X$ contains a sequence equivalent to
the ``unconditionalization'' of the basic sequence $(x_n)$. This
follows immediately from the following exact formula for the norm of
a finite sum of elementary tensors in
$\ell_1\hat\otimes_{\varepsilon} X$:
$$
\|\sum_{j=1}^n e^*_j\otimes x_j\|_\varepsilon = \sup\|\sum_{j=1}^n \pm x_j\|,
$$
where the supremum is over all choices of signs.

In the case of $\XK$ the space of compact operators $\mathcal K(\XK)$ is
isomorphic to $\ell_1\hat\otimes_\varepsilon \XK$ and so contains the
unconditionalization of any basic sequence in $\XK$. An interesting special
case is that of the basis $(d_\gamma)$; we have chosen to prove the following
proposition in a way that does not depend on the general theory of tensor
products.

\begin{prop}
The family $(e^*_\gamma\otimes d_\gamma)_{\gamma\in \Gamma}$ is an
unconditional basis of a reflexive subspace of $\mathcal K(\XK)$.
\end{prop}
\begin{proof}
Let us write $U_\gamma = e^*_\gamma\otimes d_\gamma$ considered as the rank--1
operator
$$
U_\gamma: \XK\to \XK; x\mapsto x(\gamma)d_\gamma.
$$
For a finite linear combination $W=\sum_{\gamma\in \Gamma_n } w(\gamma)
U(\gamma)$ and any  $x\in \ball \XK$ we have
$$
\|W(x) \| =\|\sum_{\gamma\in \Gamma_n} (w\gamma)x(\gamma) d_\gamma\| \le
\max_\pm \|\sum_{\gamma\in \Gamma_n} \pm w(\gamma)d_\gamma\|.
$$
We shall write $\||W\||$ for the last expression on the line above. We have
thus shown that $\|W\|\le \||W\||$.

On the other hand, if we choose $u(\gamma)=\pm 1$ for $\gamma\in \Gamma_n$  in
such a way as to achieve the maximum in the definition of $\||W\||$ and then
set $y=i_n(u)$ we have
$$
\||W\|| = \|\sum_{\gamma\in \Gamma_n} u(\gamma)d_\gamma\|=\|W(y)\|\le
\|W\|\,\|i_n\|\le 2\|W\|.
$$
Thus the operator norm $\|\cdot\|$ and the  unconditionalized norm
$\||\cdot\||$ are equivalent on $[U_\gamma:\gamma\in \Gamma]$.  It will be
convenient to work with the latter norm.

Given a linear combination $V= \sum_{\gamma} v(\gamma) U_{\gamma}$, any vector
$\sum_\gamma \pm v(\gamma)d_\gamma$ in $\XK$, (whether or not the signs achieve
the supremum in the definition of the unconditionalized norm), will be called a
{\em realization} of $W$.

If the subspace $[U_\gamma:\gamma\in \Gamma]$ is not reflexive then by
unconditionality there is a skipped block sequence equivalent to the unit
vector basis of either $c_0$ or $\ell_1$.  We shall treat the case of $\ell_1$,
leaving the (very easy) other case to the reader.

We consider a normalized skipped block sequence with $V_i=\sum _{\gamma\in
\Gamma_{p_i-1}\setminus\Gamma_{p_{i-1}}} v(\gamma)U_\gamma$ and suppose,  if
possible, that $(V_i)$ is $C$-equivalent to the usual $\ell_1$-basis for the
norm $\||\cdot\||$. More precisely, let us suppose that $\||V_i\||\le C$ for
all $i$ and that
 $$
 \||\sum_i a(i) V_i\||\ge \sum_i |a(i)|
 $$
 for all scalars $a_i$.  Let us note that if $W$ is a linear
combination of the form
 $$
W= n^{-1} \sum_{i=l+1}^{l+n} V_i,
 $$
then any realization $\hat W$ of $W$ is a $C$-$\ell_1$-average as in
Definition~\ref{ellOneDef}).  Indeed $\hat W$ is expressible as $ n^{-1}
\sum_{i=l+1}^{l+n} \hat V_i$ where the $\hat V_i$ are realizations of $V_i$,
and so satisfy $\|\hat  W_i\|\le \||W_i\||\le C$ for all $i$.

We now look at Lemma~\ref{AveRIS}.  It should be clear that, by choosing
sequences $(j_k)_{j\in \nats}$ and $(l_k)_{j\in \nats}$ growing sufficiently
fast,  we may define
 $$
W_k= n_{j_k}^{-1} \sum_{i=l_j+1}^{l_j+n_{j_k}} V_i,
 $$
in such a way that any realizations $\hat W_k$ form a $2C$-RIS in $\XK$.  In
particular
 \begin{align*}
 \||n_{j_0}^{-1}\sum_{k=1}^{n_{j_0}} W_k\|| &= \|n_{j_0}^{-1}\sum_{k=1}^{n_{j_0}} \hat
 W_k\|\\
 \intertext{for suitable realizations $\hat W_k$, yielding}
\||n_{j_0}^{-1}\sum_{k=1}^{n_{j_0}} W_k\|| &\le 12Cm_{j_0}^{-1},
\end{align*}
by Proposition~\ref{RISUpperEst}.  On the other hand,
 \begin{align*}
 \||n_{j_0}^{-1}\sum_{k=1}^{n_{j_0}} W_k\|| &=
 \||n_{j_0}^{-1}\sum_{k=1}^{n_{j_0}} n_{j_k}^{-1}\sum_{i=l_k+1}^{l_k+n_{n_k}}
 V_i\||,
\end{align*}
which is at least 1, by our assumption on $(V_i)$.

So we have a contradiction for suitably large values of $j_0$.
\end{proof}

\subsection{Open problems}

Our constructions give no clue as to whether there exists a
reflexive  Banach space on which all operators are
scalar--plus--compact.  The construction of such a space, if one
exists, will need new ideas.  We thus have no example of a reflexive
space on which all operators have non-trivial proper invariant
subspaces.  It is piquant to observe that, at the other end of the
spectrum, the construction of a reflexive space on which some
operator has no non-trivial proper invariant subspace has also
proved to be very resistant to attack.  We refer the reader to the
papers of Enflo \cite{EnOne, EnTwo} and Read \cite{ReadOne, ReadTwo}
for more about the Invariant Subspace Problem, noting the more
recent paper \cite{ReadThree} of Read, in which a strictly singular
operator is constructed which has no non-trivial proper invariant
subspace.

As we remarked in the introduction, we do not know whether an
isomorphic predual of $\ell_1$ which has the ``few-operators''
property in the scalar--plus--strictly-singular sense necessarily
also has this property in the scalar--plus--compact sense. An answer
to this would follow from an affirmative solution to the following
more general problem.

\begin{prob}
Let $X$ be a $\mathscr L_\infty$-space.  Is every strictly singular operator on
$X$ weakly compact?
\end{prob}

Working with Ch. Raikoftsalis, the present authors have recently
constructed another counterexample to the scalar--plus--compact
problem.  Like the space presented here, it is a $\mathscr
L_\infty$-space constructed by the Bourgain--Delbaen method.
However, the new space has non-separable dual and has a subspace
isomorphic to $\ell_1$.  We believe it to be the first example of an
indecomposable space containing $\ell_1$.  The only obvious
obstruction to embeddability of a given Banach space $X$ into an
indecomposable space is the existence in $X$ of a subspace
isomorphic to $c_0$.  We therefore are led to pose another problem.

\begin{prob}
Let $X$ be a separable Banach space with no subspace isomorphic to
$c_0$. Does $X$ necessarily embed in an indecomposable space?  Does
$X$ necessarily embed in a $\mathscr L_\infty$-space with the
scalar--plus--compact property?
\end{prob}

It is tempting to believe that some combination of the techniques
developed in this paper with the Bourgain--Pisier method \cite{BP}
for embedding arbitrary Banach spaces into $\mathscr
L_\infty$-spaces might provide an answer.

We should like to draw the reader's attention to the problems posed
by Bourgain \cite[page 46]{B} about the spaces $X_{a,b}$ and
$\mathscr L_\infty$-spaces in general. Problems 1,2 and 3 remain
open. Hoping that his or her appetite has been whetted by the
present paper, we leave it to the reader to find out what these
problems are. Concerning Problem 4, we now know \cite{G1} that there
is an infinite-dimensional Banach space with separable dual, no
reflexive subspace and no subspace isomorphic to $c_0$. The present
paper yields an example of an $\mathscr L_\infty$ space with no
unconditional basis sequence.  But we still do not have an example
of a space $X$ with $X^*$ isomorphic to $\ell_1$ and not containing
$c_0$ or a reflexive subspace.  Only Problem 5 has been completely
settled: each $X_{a,b}$ is saturated with $\ell_p$ for some $p$
\cite{H}.

\end{document}